\documentclass{article}
\usepackage[latin1]{inputenc}
\usepackage[T1]{fontenc}
\usepackage{fullpage}
\usepackage{amsmath}
\usepackage{amsfonts}
\usepackage{subfigure}
\usepackage{subfigmat}
\usepackage{algorithm,algorithmic}
\usepackage{amsthm}\usepackage{amssymb} 
\usepackage{url}
\newcommand{\R}{\mathbb R}
\newcommand{\Z}{\mathbb Z}

  \usepackage{graphicx}
\newcommand{\dpar}[2]{\dfrac{\partial #1}{\partial #2}}

\newcommand{\E}{{\mathcal{ E}}}

\newtheorem{theoreme}{Th\'eor\`eme}[section]

\newtheorem{remark}[theoreme]{Remark}

\usepackage[usenames,dvipsnames]{color}
\def\beq{\begin{equation}}
\def\eeq{\end{equation}}
\def\esplit{\end{split}}
\def\beqalign{\begin{array}{rl}}
\def\eeqalign{\end{array}}
\DeclareMathOperator*{\argmin}{\arg\!\min}

\makeatletter
\newcommand{\ALOOP}[1]{\ALC@it\algorithmicloop\ #1%
  \begin{ALC@loop}}
\newcommand{\ENDALOOP}{\end{ALC@loop}\ALC@it\algorithmicendloop}

\makeatother

\def\Abold{\mathbf{A}}

\def\Dbold{\mathbf{D}}
\def\Ebold{\mathbf{E}}

\def\Jbold{\mathbf{J}}

\def\Mbold{\mathbf{M}}

\def\Qbold{\mathbf{Q}}
\def\Rbold{\mathbf{R}}
\def\Sbold{\mathbf{S}}

\def\Wbold{\mathbf{W}}

\def\Zbold{\mathbf{Z}}

\def\bbold{\mathbf{b}}

\def\fbold{\mathbf{f}}

\def\mbold{\mathbf{m}}

\def\qbold{\mathbf{q}}
\def\rbold{\mathbf{r}}
\def\sbold{\mathbf{s}}
\def\tbold{\mathbf{t}}
\def\ubold{\mathbf{u}}
\def\vbold{\mathbf{v}}
\def\wbold{\mathbf{w}}
\def\xbold{\mathbf{x}}
\def\ybold{\mathbf{y}}
\def\zbold{\mathbf{z}}
\def\alphabold{\boldsymbol{\alpha}}
\def\betabold{\boldsymbol{\beta}}

\def\mubold{\boldsymbol{\mu}}
\def\rhobold{\boldsymbol{\rho}}
\def\phibold{\boldsymbol{\phi}}

\def\0bold{\boldsymbol{0}}
\def\0{\mathbf{\0}}

\newcommand{\remi}[1]{{ \color{black}{#1}}}

\begin{document}
\title{Robust Model Reduction by $L^1$-norm Minimization and Approximation via Dictionaries: Application  to Linear and Nonlinear Hyperbolic Problems}
\author{Remi Abgrall $^{(*)}$ and David Amsallem$^{(\dagger)}$ \\
(*) Institut f\"ur Mathematik, Winterthurstrasse 190, \\CH 8057  Z\"urich, Switzerland\\
($\dagger$) Department of Aeronautics and Astronautics, 
496 Lomita Mall, \\Stanford University, Stanford, CA 94305-3035, USA}
\date{\today}
\maketitle
\graphicspath{{/Users/abgrall/Reduit/Papier_David/}}
\begin{abstract}
We propose a novel model reduction approach for the approximation of non linear hyperbolic equations
in the scalar and the system cases. The approach relies on an offline computation of a dictionary of solutions together with an online $L^1$-norm minimization of the residual. It is  shown why this is a natural framework for hyperbolic problems and tested on nonlinear problems such as Burgers' equation and the one-dimensional Euler equations involving shocks and discontinuities. Efficient algorithms are presented for the computation of the $L^1$-norm minimizer, both in the cases of linear and nonlinear residuals. Results indicate that the method has the potential of being accurate when involving only very few modes,  generating physically acceptable, oscillation-free, solutions.
\end{abstract}
\section{Introduction}

Many engineering applications require the ability to simulate the behavior of a physical system in real-time. This requirement holds in particular when a full parametric exploration of the behavior of the system is sought. In aerodynamics, such an exploration can be done to compute the flow around an aircraft for varying boundary conditions or to design its shape  to maximize lift and minimize drag. Uncertainty quantification also requires a large number of simulations with varying parameters in order to propagate chaos by means of a Monte-Carlo method or calibrating input parameters by a Markov chain technique. A third important application is flow control.

When such a large number of simulations is required, the cost of one simulation is critical to the application at hand. This cost can be lowered by using sophisticated computer science techniques such as parallelization but such techniques are usually not enough to allow full parametric exploration, especially when computational resources are limited. 

Alternatively, model reduction techniques can alleviate the cost of such repeated simulations with limited computational resources~\cite{legresley00,buithanh08,amsallem10,amsallem14:smo}.  Model reduction is directly based on the underlying high-dimensional model (HDM) that results from a standard finite element, finite volume of finite differences formulation. In the present paper, 
Partial Differential Equations (PDE) of the following type are considered:
\begin{equation}
\begin{split}
\label{eq:1}\dpar{U}{t}+L(U)=0 & \qquad x\in \Omega,~t\in [0,T]\\
B(U)=g & \qquad x\in \partial \Omega,~t\in [0,T]\\
U(x,t=0)=U_0(x) & \qquad x\in \Omega
\end{split}
\end{equation}
$L$ is a differential operator (for example the Laplacian or the divergence of a flux), and $B$ a boundary operator. In this paper, we are particularly interested in the
case where the solution $U(x,t)\in \R^p$ is a scalar or a vector and $L$ is the divergence of a flux $F$.
Two examples will be considered by increasing order of complexity:
\begin{itemize}
\item  Burgers' equation for which $U=u$ is scalar:
\begin{itemize}
\item Its unsteady version,
$$\dpar{u}{t}+\dpar{}{x}\left(\frac{1}{2}u^2\right)=0, \qquad u(x,0)=u_0(x)$$
with periodic boundary conditions
\item It steady version  with weak Dirichlet boundary conditions
\end{itemize}
\item The one-dimensional compressible Euler equations for which $U=(\rho, \rho u, E)$, $F(U)=(\rho u, \rho u^2+p, u(E+p))$ and  the perfect gas equation of state holds:
$$p=(\gamma-1)\left( E-\frac{1}{2}\rho u^2\right).$$
$\rho$ denotes the density, $u$ the velocity, $p$ the pressure and $E$ the energy. 
\end{itemize}
After discretization in space,  the solution is denoted as $\ubold(t)\in\mathbb{R}^{Np}$. The PDE is here parameterized by a parameter vector $\mubold\in\mathbb{R}^m$ that allows changes in the operator $L$, the boundary operator $B$ or the initial conditions. For simplicity and without loss of generality, this parametric dependency will be omitted in the next paragraphs.

Instead of allowing any value of the solution degrees of freedom $\ubold$, model reduction however restricts the solution to be contained in a subspace of the underlying high-dimensional space. This subspace is determined by an optimized reduced basis that is determined in a training phase. Thus, a large number of degrees of freedom (say millions) are represented by only a few number of coefficients in the representation of the full solution in terms of the reduced basis vectors, leading to important computational savings. Two important questions arise at this point: (1) How can an optimal reduced basis be constructed? and (2) How can the evolution of the reduced coefficients be computed in a stable fashion?

A popular method for choosing an ``optimal'' basis is Proper Orthogonal Decomposition (POD), first introduced as a tool for the analysis of flows by
 Lumley~\cite{Lumley} and then extended and popularized by Sirovich~\cite{sirovich87}. The idea behind POD is to collect a few snapshots of the solution and then compute the best approximation of these snapshots in terms of a small number of reduced basis vectors. Mathematically speaking, if $\ubold_i(t_l)\in\mathbb{R}^p$ denotes the value of the discrete solution $\ubold$ at grid point $\mathbf{x}_i,~i=1,\cdots,N$ and at time $t_l,~l=1,\cdots,N_t$, POD constructs $M$ orthogonal functions $\phibold_\ell \in \big [L^2(\R^d)\big]^p$ such that the following functional  is minimized:
 \begin{equation}
 \mathcal{J}(\phibold_1,\cdots,\phibold_M) = \sum_{l=1}^{N_t}\sum_{i=1}^{Np}\left\|\ubold_i(t_l)-\sum_{\ell=1}^M \langle u(t_l),\mathbf{\phi}_{\ell}\rangle \mathbf{\phibold}_{\ell i}\right\|_2^2,
 \end{equation}
 where ${\phibold}_{\ell i}\in\mathbb{R}^p$ denotes the value of $\phibold$ at $\mathbf{x}_i$. $\|~\cdot~\|$ denotes here the Euclidean norm in $\R^p$, and $\langle ~\cdot~,~\cdot ~\rangle $ is the $L^2$ norm. A minimum of the functional $\mathcal{J}$ can be analytically computed by Singular Value Decomposition \cite{Golub}, and the reduced basis vectors $\phibold_{\ell}$ are the left singular vectors of the snapshots matrix
$$\Sbold=\begin{pmatrix}
\mathbf{u}_1(t_1) & \ldots & \mathbf{u}_1(t_{N_t})\\
\vdots & \vdots & \vdots \\
\mathbf{u}_N(t_1) & \ldots & \mathbf{u}_N(t_{N_t})
\end{pmatrix}.$$ 
Defining by $\{\lambda_\ell\}_{l=1}^{N_t}$ the positive eigenvalues of $\Sbold^T\Sbold$ sorted decreasingly, the error associated with the minimum of the functional is
\begin{equation}
\mathcal{J}(\phibold_1,\cdots,\phibold_M)  = \sum_{\ell=M+1}^{N_t} \lambda_\ell.
\end{equation}

In the continuous case, the functions $\mathbf{\phibold}_{\ell}(\xbold)\in\mathbb{R}^p$, 
 are the solution of Fredholm alternative
\begin{equation}\int_{\Omega} R(\mathbf{x},\mathbf{x}')\mathbf{\phibold}_\ell(\mathbf{x}')d\mathbf{x}'=\lambda_\ell \mathbf{\phibold}_\ell(\mathbf{x}), \qquad \text{for all }\mathbf{x}\in \Omega,
\label{evalue}\end{equation}
where $\Rbold(\xbold,\xbold') = \ubold(\xbold)\ubold(\xbold')^T$.

In both the discrete and continuous cases, the basis dimension $M$ is determined on the basis of the decay of the eigenvalues $\lambda_\ell$. Given a tolerance $\epsilon\ll 1$, $M$ is selected as the smallest dimension such that the following relative truncation error is smaller than $\epsilon$.
\begin{equation}
\frac{\mathcal{J}(\phibold_1,\cdots,\phibold_M) }{\sum_{l=1}^{N_t}\sum_{i=1}^{Np}\left\|\ubold_i(t_l)\right\|_2^2} = \frac{\sum_{\ell=M+1}^{N_t} \lambda_\ell}{\sum_{\ell=1}^{N_t} \lambda_\ell}.
\end{equation}

\begin{figure}
\begin{center}
{\includegraphics[width=0.45\textwidth,clip=]{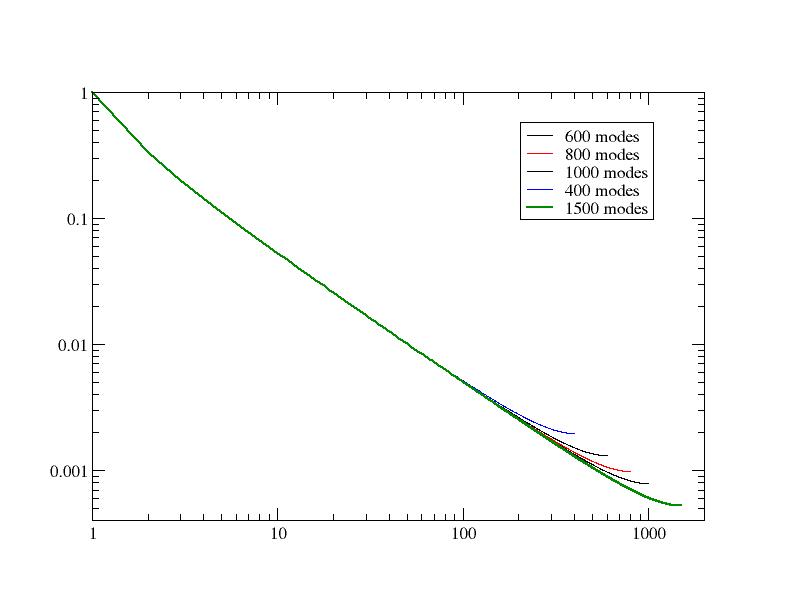}}
\end{center}
\caption{\label{fig:evalue}\emph{
Ratio of POD eigenvalues $\log(\lambda_k(N)/\lambda_1(N))$ for $N=400, 600, 800, 1000, 1500$ grid points.}}
\end{figure}

 In general, one expects the eigenvalues $\lambda_\ell$ to decrease very rapidly to $0$. This allows, when this assumption is true, to consider only the most energetic modes in the decomposition. Unfortunately, it is not always the case that  the eigenvalues $\lambda_\ell$ are rapidly converging to zero. This is demonstrated by the following simple counter example for which a simple scalar advection problem defined on $\Omega = [0,1[$ is considered:
\begin{subequations}\label{pde}
\begin{equation}
\label{pde:1}\dpar{u}{t} +\dpar{u}{x}=0
\end{equation}
with the boundary condition 
\begin{equation}
\label{pde:2}u(0,t)=1
\end{equation}
and the initial condition 
\begin{equation}
\label{pde:3}u(x,0)=0.
\end{equation}
 The solution is given by a traveling discontinuity
$$u(x,t)=\left\{\begin{array}{ll}
1 & \text{if } x\leq \min(t,1)\\
0& \text{otherwise.}
\end{array}\right .$$
\end{subequations}
Considering grids $x_i=i/N$, $i=0, \ldots, N$ for varying number of grid points $N$ and  snapshots collected at times as $t_k=k\Delta t$, with $\Delta t=1/N$, a series of POD bases is constructed numerically. For each grid size $N$, the eigenvalues $\lambda_\ell(N)$ are reported in Figure~\ref{fig:evalue}. One can observe that the ratio $\lambda_\ell(N)/\lambda_1(N))$  behaves like $1/k$ and $\max(\lambda_\ell)$ behaves like $0.63 \, N$. This illustrates that it is not possible to select only a few dominant modes, due to the slow decay of the POD eigenvalues.  This example also illustrates why most of the work on model 
reduction has been focused on regular problems, and for fluids, on incompressible flows, see e.g. among many others \cite{zbMATH01784252,zbMATH05066006,Veroy}.
For compressible (but regular) flows, one of the early work is \cite{Rowley2004}, then one may mention \cite{GNAT} for compressible turbulent flows, \cite{amsallem14:morepas} for compressible inviscid flows and~\cite{amsallem08,Barone2008,Kalashnikova2009,amsallem10} for the case of linearized compressible inviscid flows.

%

Concerning compressible fluids, there is another difficulty. In problem~\eqref{evalue}, one needs a norm. In the case of incompressible flows, a natural norm is related
to the kinetic energy. For compressible materials,  however, one needs to take into account the density, velocity and the energy, i.e. the thermodynamics. A simple $L^2$-norm cannot be
 used because one cannot combine in a quadratic manner these variables, for dimensional reasons. Only a non-dimensionalization of the variables~\cite{lesoinne01} can alleviate the dimensionality issue~\cite{amsallem08,GNAT,amsallem14:morepas}.
 
 The natural equivalent of the $L^2$-norm is however related to the entropy, which 
is not quadratic: if a minimization problem can be set up, its solution is non trivial. These arguments were raised in \cite{Rowley2004}, and an energy-based norm was developed  in~\cite{Barone2008,Kalashnikova2009} for linearized compressible flows. 

To circumvent those issues, an approach based on a dictionary of solutions~\cite{maday12,kaulmann13} is developed in this work as an alternative to using a truncated reduced basis based on POD.
The elements of this dictionary are solutions $\ubold(t_l;\mubold_j)$ computed for varying values of time $t_l$ and parameter $\mubold_j\in\mathbb{R}^m$. Selecting appropriate parameter samples $\mubold_j\in\mathcal{D} \subset \mathbb{R}^m$ is a crucial step that can affect the accuracy of the reduced-order model in the parameter domain. Greedy sampling procedures have been developed when error estimates are known~\cite{Prudhomme,Veroy,zbMATH02136393,zbMATH06198352,zbMATH05486455,zbMATH01784252,pdt14}. In this work, the issue of optimal sampling is not considered the main focus of the paper lies in establishing an effective model reduction approach based on dictionaries for hyperbolic problems.
 The development of a strategy to sample optimal values of $\mubold$ in this context will be the topic of further research.

In addition to choosing an appropriate dictionary $\mathcal{D}$, selecting an approach for computing a reduced solution based on that dictionary is also crucial. For self-adjoint systems, Galerkin projection is a natural approach but there is no motivation for using Galerkin projection for nonlinear compressible flows. Instead, strategies based on the minimization of the residual arising from the reduced approximation have been successfully developed for compressible flows in~\cite{legresley00,buithanh08,GNAT,amsallem14:morepas}. These approaches rely on a minimization of the residual in the $L^2$ sense. In the present work, this minimization problem is extended to the more general minimization using a $L^q$-norm, with emphasis on $q=1$ and $q=2$. For nonlinear systems, an additional step, hyper-reduction, is required to ensure an efficient solution of the reduced system~\cite{Ryckelynck2005,GNAT,amsallem14:morepas}. Hyper-reduction is not considered in this work but will be the subject of follow-up work.

This paper is organized as follows. Motivations for using the $L^1$-norm in the case of hyperbolic systems are given in Section~\ref{sec:L1} \remi{where we show that $q=1$ is very closely linked to the concept of weak solutions of  hyperbolic problems. } The proposed model reduction approach is then developed in Section~\ref{formulation} in both the steady and unsteady cases. Finally, the proposed procedure is applied  to the model reduction of several steady and unsteady systems in Section~\ref{sec:appli} and conclusions are given in Section~\ref{sec:conclu}.

\section{Motivation for the $L^1$-norm.}\label{sec:L1}
\remi{In solving minimization problems, it is quite usual to minimize a residual with respect to the $L^q$-norm for a suitable $q$. The choice $q=2$ is very common because it amounts to minimize in some least-squares sense and many efficient algorithms are available. In the case of hyperbolic problems, as we are concerned with here, this is still a convenient choice (after proper non-dimensionalization as mentioned above), but it might not be the most natural one, as demonstrated in the work of Guermond et al. on Hamilton Jacobi equations and transport problems \cite{Guermond,guermond2}. In particular these works show, at least experimentally, that the numerical solution has an excellent non-oscillatory behavior by minimizing the $L^1$-norm of the PDE residual. In fact, this observation is our original motivation for choosing the $L^1$-norm, since we are interested in preserving the non-oscillatory nature of solutions. In this section, we further justify the choice of the $L^1$-norm applied to the residual, and show that it is closely related to the weak formulation of the problem.}
The following discussion is formal. 

Let us consider the problem
\begin{equation}
\label{eq:2}
\dpar{U}{t}+\text{ div } {F}(U)=0
\end{equation}
defined on $\Omega \subset \R^d$, $t> 0$. The steady problem can be done in the same exact manner. We assume that the solution $U$ belongs to $\R^p$, so that $F=(F_1, \ldots, F_p)^T$.
The weak form of this is: for any $\varphi \in \left[C^1(\Omega)\right]^p$ and with compact support, we have:
$$\int_\Omega \varphi(x,t)\bigg ( \dpar{U}{t}+\text{ div } F(U) \bigg )dx=0.$$
Integrating by parts yields
$$
\int_\Omega \dpar{\varphi}{t} U dx+\int_\Omega \nabla \varphi \cdot F(U) dx=0.$$
If we restrict ourself to the set of test functions $\left\{\varphi\in \left[C^1(\Omega)\right]^p, ||\varphi||_{\infty}\leq 1\right\}$, we have that $U$ is a solution if:
$$\sup\limits_{\{\varphi\in \left[C^1(\Omega)\right]^p, ||\varphi||_{\infty}\leq 1\}}\Bigg (\int_\Omega \dpar{\varphi}{t} U dx+\int_\Omega \nabla \varphi \cdot F(U) dx\Bigg )=0.$$
Let us now real the definition of the total variation of a function $g\in L^1(\R^d)$:
$$TV(g)=\sup\limits_{\varphi\in C^1_0(\R^d)\cap L^\infty(\R^d), ||\varphi||_{\infty}\leq 1} \left\{\int_{\R^d} \nabla \varphi (x)\cdot g(x) dx\right\},$$
and we see that if in addition $g\in C^1(\R^d)$, $TV(g)=\int_{\R^d}||\nabla g||dx=||\nabla g||_{L^1(\R^d)}.$

Before going further, let us mention the following classical result that  will be useful.
Consider $\{x_i\}_{i\in \Z}$ a strictly increasing sequence in $\R$, we define $x_{i+1/2}=\frac{x_i+x_{i+1}}{2}$. We assume that $\R=\cup_{i\in \Z}[x_{i-1/2},x_{i+1/2}[$ and consider $g$ defined by: for any $i\in \Z$, 
$$g(x)=g_i \text{ if } x\in [x_{i-1/2},x_{i+1/2}[,$$
we see that
$$TV(g)=\sum_{i\in \Z} |g_{i+1}-g_i|.$$

Thanks to this definition, we see that
if we define the space-time flux $\mathcal{F}=(U,F)$, $U$ is a weak solution if and only if the total variation of $\mathcal{F}$ vanishes,
$ TV\big ( \mathcal{F}) =0.$

\medskip

Now, instead of having the exact solution, we consider an approximation procedure that enables, from $\ubold^{n}\approx U(~.~, t_n)$, to compute $\ubold^{n+1}\approx U(~.~, t_{n+1})$,
say $\mathcal{L}(\ubold^n, \ubold^{n+1})$.

For instance, assume that  we have a finite volume method and $d=1$: for any grid point $i\in\{1,\cdots,N\}$,
$$
\left[\mathcal{L}(\ubold^n, \ubold^{n+1})\right]_i= \Delta x(\ubold_i^{n+1}-\ubold_i^n)+ {\Delta t}\big (\fbold_{i+1/2}(\ubold^n)-\fbold_{i-1/2}(\ubold^n)\big ).$$

A way to evaluate $\ubold^{n+1}$ is to minimize the total variation, i.e.
$$TV(\mathcal{L})=\sum_{i\in \mathcal{I}} \Big |\Delta x (\ubold_i^{n+1}-\ubold_i^n)+ {\Delta t}\big(\fbold_{i+1/2}(\ubold^n)-\fbold_{i-1/2}(\ubold^n)\big )\Big |,$$
$$\ubold^{n+1}= \argmin _{\vbold \text{ piecewise constant }} \sum_{i\in \mathcal{I}} \bigg |\Delta x (\vbold_i-\ubold_i^n)+ \Delta t(\fbold_{i+1/2}(\ubold^n)-\fbold_{i-1/2}(\ubold^n))\bigg |.$$
Clearly, if $\mathcal{I}$ is equal to the set of grid points, the solution is given by 
$$\ubold_i^{n+1}=\ubold_i^n -\dfrac{\Delta t}{\Delta x} \bigg (\fbold_{i+1/2}(\ubold^n)-\fbold_{i-1/2}(\ubold^n)\bigg )$$
and nothing new is gained. 
 When $\mathcal{I}$ is not equal to the set of degrees of freedom, then something new happens. We expect precisely to exploit this idea, or ideas related to this.

In the remainder of this paper, this idea is exploited in the case of model reduction, for which $\mathcal{I}$ is not equal to the set of grid points and the $TV$ semi-norm slightly modified in order to guarantee (1) that a unique solution to the minimization problem exists, and (2) that the minimization problem is as easy as possible to solve.

\section{Formulation}\label{formulation}

\subsection{High-dimensional model}

Without loss of generality, the case of the  classical finite volume method is considered to define the High Dimensional Model (HDM). 
A computational domain $\Omega\subset\mathbb{R}^d$ is considered, and in most of this paper, $\Omega\subset\R$, that is $d=1$. Starting from a subdivision
$\ldots  < x_j <x_{j+1}< \dots$, we construct control volumes $K_{j}=[x_{j-1/2},x_{j+1/2}[$, $j\in \Z$ where
$$
x_{j+1/2}=\frac{x_j+x_{j+1}}{2}.
$$
A finite volume semi-discrete formulation of \eqref{eq:1} writes 
\begin{subequations}
\label{edp}
\begin{equation}\label{edp:1}
|K_j| \dfrac{d\ubold_j}{dt}+\fbold_{j+1/2}(\ubold)-\fbold_{j-1/2}(\ubold)=0
\end{equation}
where $\fbold_{j+1/2}$ is a consistent numerical flux. In each applications, we consider Roe's formulation and a first order scheme.
We assume either compactly supported initial conditions or initial conditions with periodicity
\begin{equation}
\label{edp:2}
\ubold_j(t=0)\approx \dfrac{1}{|K_j|}\int_{K_j} U_0(x)dx.
\end{equation}
In \eqref{edp:1}, $\ubold_j$ stands for an approximation of the average of the solution in the cell $K_j$,
$$\ubold_j(t)\approx\dfrac{1}{|K_j|}\int_{K_j} U(x,t) dx.
$$
\end{subequations}
The time stepping is done in a  standard way, for instant by Euler time stepping.

\subsection{Model Reduction by residual minimization over a dictionary}
\subsubsection{Steady problems}\label{sec:steady}
Two approaches are available to solve a steady state associated with problem \eqref{eq:1}. The first one is to use a homotopy approach~\cite{kelley98} with pseudo-time stepping, resulting in the solution of an unsteady problem which limit solution is the desired steady state. The procedure described for unsteady systems in Section~\ref{sec:unsteady} can be, in principle applied to this case. The second approach is by a direct solution of the steady-state problem. The discretized steady-state problem writes 
$$\rbold(\ubold(\mubold),\mubold) = 0$$ 
where $\rbold(\cdot,\cdot)$ is usually a nonlinear function of its arguments, referred to as the residual. This set of nonlinear equations is typically solved by Newton-Raphson's method. This second approach is followed in this work for steady problems.

The parameter vector $\mubold\in\mathcal{P}\subset\mathbb{R}^m$ can, for instance, parameterize the boundary conditions associated with the steady-state problem. The parametric domain of interest $\mathcal{P}$ is assumed here to be a bounded set of $\mathbb{R}^m$.

The solution manifold $\mathcal{M} = \left\{\ubold(\mubold)~\text{s.t}~ \mubold\in\mathcal{P}\subset\mathbb{R}^m\right\}$ is assumed to be of small dimension. This manifold $\mathcal{M}$ belongs to $L^\infty(\R^d)\cap BV(\R^d)$, and thus can be locally described by some mapping $\theta: \mathcal{P}\mapsto L^\infty(\R^d)\cap BV(\R^d)$.
To approximate this mapping, we consider a family of $r$ parameters in $\mathcal{P}$, $\{\mubold_\ell\}_{\ell=1}^r$, and compute the associated dictionary of solutions  $\mathcal{D} = \left\{\ubold(\mubold_\ell)\right\}_{\ell=1}^r$ of~\eqref{edp}. 

The steady-state $\ubold(\mubold)$ is then approximated as a linear combination of the pre-computed dictionary elements $\mathcal{D}$ as
\begin{equation}
\ubold(\mubold) \approx \sum_{\ell=1}^r \alpha_\ell(\mubold) \ubold(\mubold_l).
\end{equation}

For a new value of the parameters $\mubold\in\mathcal{P}$, the reduced coordinates  $\left\{\alpha_\ell(\mubold)\right\}_{\ell=1}^r$ are then computed as the solution of the minimization problem
\begin{equation}
\alphabold(\mubold) :=(\alpha_1(\mubold), \ldots , \alpha_r(\mubold))= \argmin_{\betabold=(\beta_1,\cdots,\beta_r)}J\left( \rbold\left(\sum_{\ell=1}^r \beta_\ell\ubold(\mubold_l),\mubold\right),\betabold\right).
\end{equation}

In this paper we consider for $J$ the following convex functionals, which are described in more details in Appendix~\ref{algo}.
\begin{itemize}
\item the $L^2$-norm $J(\zbold,\xbold)=\|\zbold\|_2$ and its regularized version $J(\rbold,\betabold)=\|\rbold\|_2+\eta \|\betabold\|_2$ with $\eta >0$. 

When $\rbold$ is a linear function of $\betabold$ ($\rbold = \Abold\betabold+\bbold$), this choice of functional results in the solution of an ordinary least-squares problem, described in Appendix~\ref{L2:L}. When  $\rbold$ is a nonlinear function of $\betabold$, Gauss-Newton or Levenberg-Marquardt procedures can be used to minimize $J$, as described in Appendix~\ref{L2:NL}.

\item the $L^1$-norm $J(\rbold,\betabold)=\|\rbold\|_1$  or its regularized variant,  $J(\rbold,\betabold)=\|\rbold\|_1+\eta \|\betabold\|_1$ with $\eta>0$.

 Two approaches are considered  to minimize $J$ when $\rbold$ is a linear function of $\betabold$.
\begin{enumerate}
\item Linear Programming (LP), involving the solution of an optimization problem with $2m+r$ variables and $3m$ constraints.  
\item The Iteratively Reweighted Least-Squares approach (IRLS)~\cite{daubechies08}.
\end{enumerate}
Both approaches are described in great detail in Appendix~\ref{app:L1lin}. When $\rbold$ is a nonlinear function of $\betabold$, a Gauss-Newton-like procedure can be used in combination with either the LP or IRLS approaches, as described in Appendix~\ref{app:L1nonlin}.
 Unicity of the solution can be guaranteed by setting the regularization term $\eta>0$.
\item The Huber function $J(\rbold) = \sum_{i=1}^m \phi_M(r_i)$~\cite{huber11} as described in Appendix~\ref{ssec:Huber}. In the present work, minimization of the Huber functional is carried out by the IRLS approach. The procedure is described in details in Appendix~\ref{ssec:Huber}.

The Huber functional can also be regularized by defining $J(\rbold,\betabold) = \sum_{i=1}^m \phi_M(r_i) + \eta\|\betabold\|_q$ with $q=1$ or $q=2$ and $\eta>0$.
\end{itemize}

\begin{remark}\label{remark:label}
\begin{itemize}
\item Decreasing the dimensionality of the solution space from $N$ to $r$ is not enough to gain computational speedup when the system to be solved is nonlinear. An additional level of approximation, hyper-reduction, is necessary~\cite{chaturantabut10,carlberg11,GNAT,amsallem12:localROB}. Hyper-reduction is not considered in the present work, but will be the focus of future work.
\item A careful selection of the sample parameter samples $\left\{\mubold_\ell\right\}_{\ell=1}^r$ is necessary in order to generate a reduced-order model that is accurate in the entire parameter domain $\mathcal{P}$. Greedy sampling techniques~\cite{Prudhomme,Veroy,zbMATH02136393,zbMATH06198352,zbMATH05486455,zbMATH01784252,pdt14}, associated with a posteriori error estimates, have been successfully used to construct reduced models that are robust and accurate in a parameter domain $\mathcal{P}$. These techniques are not considered in this paper but will also be the focus of future work. 
\end{itemize}
\end{remark}


\subsubsection{Unsteady problems}\label{sec:unsteady}

For simplicity, in the remainder of this section, we assume that only the initial condition $\ubold^0(\mubold)$ depends on a parameter vector $\mubold\in \mathcal{P}\subset\mathbb{R}^m$. Again, the family of solutions $\ubold(\mubold)$ of the Cauchy problem \eqref{edp} is then conjectured to belong to a low dimensional manifold $\mathcal{M}$ when the initial condition is parameterized in \eqref{edp:2}. 

To approximate this mapping, we consider a family of $r$ parameters in $\mathcal{P}$, $\{\mubold_\ell\}_{\ell=1}^r$, and compute the associated solutions of~\eqref{edp} for respective initial conditions $\ubold^0(\mubold_\ell)$, $\ell=1, \ldots , r$. 

Once these solutions are computed, we propose to approximate, for any parameter $\mubold\in \mathcal{P}$ the solution $\{\ubold^n(\mubold)\}_{n=0}^{N_t}$ associated with an initial condition $\ubold^0(\mubold)$ by approximating it  as
$$\ubold^n(\mubold) = \sum_{\ell=1}^r \alphabold^n(\mubold)\ubold^n(\mubold_\ell) $$
by the following procedure:
\begin{enumerate}
\item Initialization: determine the reduced coefficients $\{\alpha^0_\ell(\mubold)\}_{\ell=1}^r$ as: 
$$\alphabold^0(\mubold):=(\alpha_1^0(\mubold), \ldots , \alpha_r^0(\mubold))= \argmin_{\betabold=(\beta_1,\cdots,\beta_r)}J\bigg ( \sum_{\ell=1}^r \beta_\ell \ubold^0({\mubold_\ell}) ,\betabold\bigg),$$
for a given choice of functional $J(\ubold,\betabold)$. 
\item Assume that $\alphabold^n(\mubold) = (\alpha_1^n(\mubold), \ldots , \alpha_r^n(\mubold))$ is known, determine $\alphabold^{n+1} = (\alpha_1^{n+1}, \ldots , \alpha_r^{n+1})$ such that:
\begin{equation*}
\begin{split}{\alphabold}^{n+1}(\mubold)= \argmin_{\betabold=(\beta_1,\cdots,\beta_r)}J\Bigg (  \sum_{\ell=1}^r \beta_\ell \ubold^{n+1}({\mubold_\ell})-\wbold^n(\mubold) -\dfrac{\Delta t}{\Delta x} \bigg (\fbold_{1/2}(\wbold^n)-\fbold_{-1/2}(\wbold^n)\bigg ),\betabold\Bigg )
\end{split}
\end{equation*}
where
$$\wbold^n(\mubold)=\sum_{\ell=1}^r \alpha^n_\ell(\mubold) \ubold^{n}({\mubold_\ell}).$$
\end{enumerate}
We see that the second step can be written as: find $\alphabold^{n+1}(\mubold)$ that minimizes
$$\overline{J}( \Abold^{n+1}\alphabold^{n+1} -\bbold^n):=J(\Abold^{n+1}\alphabold^{n+1}-\bbold^n, \alphabold^{n+1})$$
where the matrix $\Abold^{n+1}$ can be written by blocks as
\begin{equation}
\Abold^{n+1} = \begin{pmatrix}
\ubold^{n+1}_1(\mubold_1) & \ldots &\ubold^{n+1}_1(\mubold_r)\\
\vdots & \vdots & \vdots \\
\ubold^{n+1}_N(\mubold_1) & \ldots & \ubold^{n+1}_N(\mubold_r)
\end{pmatrix}
\end{equation}
 and  $\bbold^n$ depends on known data.

A few immediate remarks can be made. 
\begin{remark}\label{remark:label2}
\begin{itemize}
\item In the case of a linear flux, Problem~\eqref{eq:1} is linear. If $\mathcal{S}_t$ is the mapping between the initial condition $U_0$ and the solution at time $t$, we have $\mathcal{S}_t(U+V)=\mathcal{S}_t(U)+\mathcal{S}_t(V)$. This means the exact solution of the Cauchy problem  with 
$U_0=\sum_\ell \alpha_\ell^0 U_0(\mubold_\ell)$ is $\mathcal{S}_t(U_0)=\sum_\ell \alpha_\ell^0 \mathcal{S}_t(U_0({\mubold_\ell}))$. In the case of a linear scheme, 
minimizing the functional $J$ should result in $\alphabold^n=\alphabold^0$ for any $n\geq 0$.
\item
In the case of an explicit background scheme, the choice of the numerical flux, how high order is reached, and the choice of time stepping has no influence on the overall procedure: any sub-time step would be treated similarly. In this paper, we have chosen a first order method with Euler time stepping in the case of unsteady problem. 
\item In the case of an implicit scheme, a Newton-like procedure can be applied to minimize the functional as in~\cite{GNAT}. At each time step, the procedure is then identical as in the steady case described in Section~\ref{sec:steady}.
 \end{itemize}
\end{remark}


\section{Numerical examples}\label{sec:appli}
\subsection{{Simple regression}}\label{ssec:regression}
{As a first example, the choices of functionals proposed in Section~\ref{sec:steady} are applied to a very simple regression problem. This example illustrates the behavior of each approach. In a first case, $22$ points $\{x_i\}_{i=1}^{22}$ are randomly selected in the interval $[0,1]$ and the coordinates $\{y_i\}_{i=1}^{22}$ are drawn from a distribution $2x_i+0.4+ 10^{-1}\; \mathcal{U}(-1,1)$ where $\mathcal{U}(-1,1)$ denotes the uniform distribution between $-1$ and $1$. The regression approximation is therefore $y\approx \alpha_1x+\alpha_2$ and the target is $\alphabold^\star = (\alpha^\star _1,\alpha^\star _2)=(2,0.4)$.

The  $L^1$-norm minimization by LP and IRLS, $L^2$-norm and Huber function minimizations are then used as functionals for that problem. The results are reported in  Figure~\ref{fig:Regression}(a) and Table~\ref{tab:RegressionNoOutliers}. One can observe that all four methods provide a good approximation of $\alphabold^\star$. Furthermore, the two $L^1$ minimization procedures as well as the Huber function minimization return identical results.

\begin{table}[htdp]
\begin{center}\begin{tabular}{|c|c|c|c|c|c|}
\hline
 & Target & $L^2$-norm & $L^1$-norm (LP) & $L^1$-norm (IRLS) & Huber function \\
 \hline
 $\alpha_1$ &2 & 1.9520 & 1.9037 & 1.9037 & 1.9037 \\
 \hline
  $\alpha_2$  &0.4 & 0.4087 & 0.4408 & 0.4408 & 0.4408\\
  \hline
 \end{tabular} 
 \caption{Regression: results without outliers}
 \label{tab:RegressionNoOutliers}
\end{center}
\end{table}
 
 \begin{figure}
\begin{center}
 \begin{subfigmatrix}{2}
  \subfigure[Case without outliers]{ \includegraphics[width=4.25in]{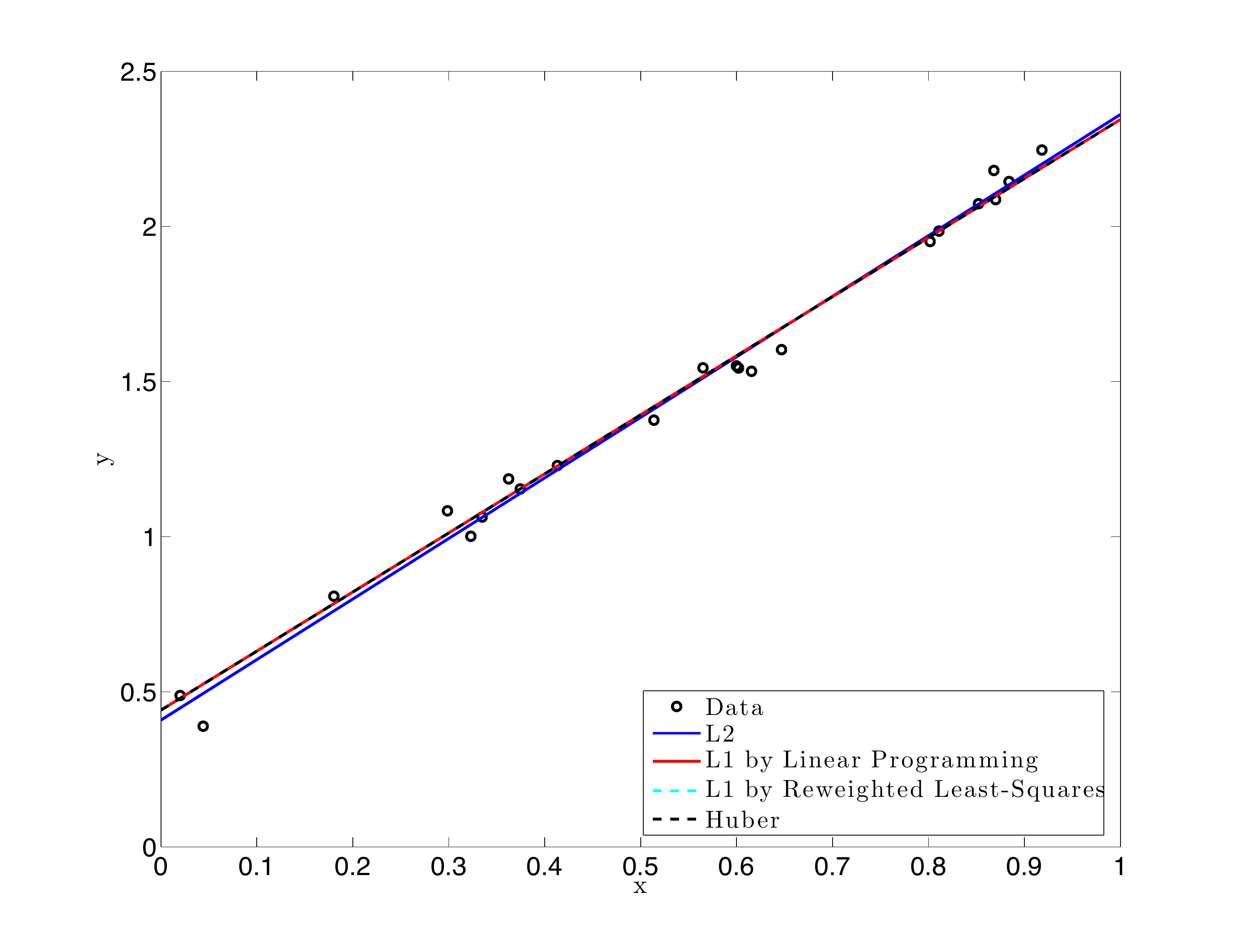}}
  \subfigure[Case with outliers]{ \includegraphics[width=4.25in]{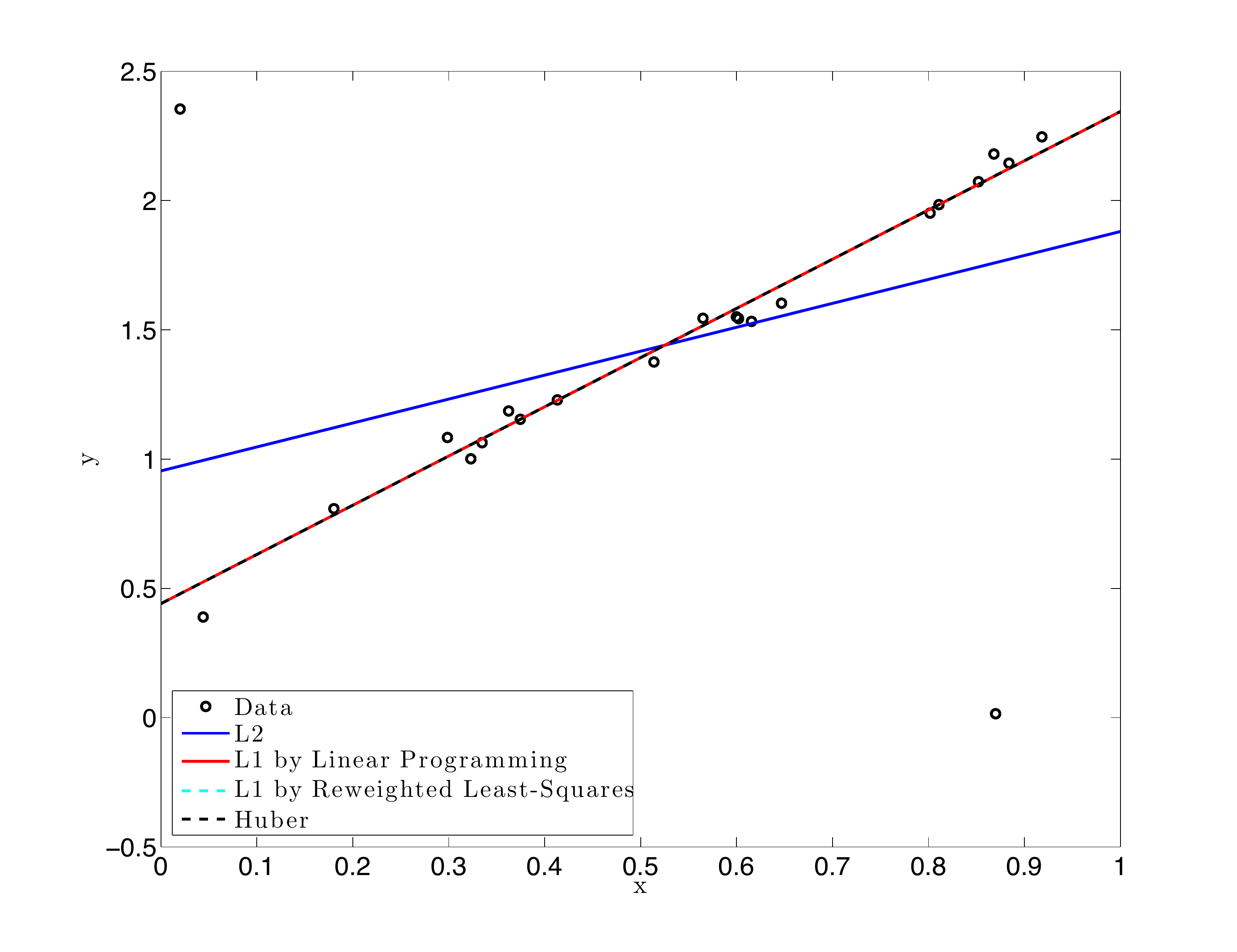}}
  \end{subfigmatrix}
 \caption{Regression: Comparison of the methods}
 \label{fig:Regression}
 \end{center}
\end{figure}

 \begin{table}[htdp]
\begin{center}\begin{tabular}{|c|c|c|c|c|c|}
\hline
 & Target & $L^2$-norm & $L^1$-norm (LP) & $L^1$-norm (IRLS) & Huber function \\
 \hline
 $\alpha_1$ &2 & 0.9256 & 1.9037 & 1.9037 & 1.9037 \\
 \hline
  $\alpha_2$  &0.4 &0.9545 & 0.4408 & 0.4408 & 0.4408\\
  \hline
 \end{tabular} 
 \caption{Regression: results with outliers}
 \label{tab:RegressionOutliers}
\end{center}
\end{table}
 
 In a following set of experiments two outliers points are defined and all four approaches applied to that new regression problem. The results are reported in  Figure~\ref{fig:Regression}(b) and Table~\ref{tab:RegressionOutliers}. One can observe that the $L^2$-norm minimization  procedure is much more sensitive to the outliers. As such, the regression coefficients returned by that procedure differ greatly from the previous case and $\alphabold^\star$ is inaccurately estimated. On the other hand, the $L^1$-norm and Huber minimization procedures are much less sensitive to the outlier points and accurate estimations of $\alphabold^\star$ are provided. Again, the three estimations are identical.
}


\subsection{{Model reduction of steady problems}}

\subsubsection{{One-dimensional advection equation}}
{ The following one-dimensional steady advection equation is considered:
\begin{equation}
\frac{d{u}}{dx}(x) = f(x;\mu),~x\in[0,~1], 
\end{equation}
where $$f(x;\mu) = \frac{-2k\exp(-2k(x-\mu))}{(1+\exp(-2k(x-\mu)))^2}$$ and $k=100$. The solution exhibits a sharp gradient at location $x=\mu$ similar to a shock. A Dirichlet boundary condition $u(0) = 1$ is imposed at $x=0$. This PDE is discretized by finite differences using a uniform mesh, resulting in a HDM of dimension $N=10^3$.

\begin{figure}
\begin{center}
{\includegraphics[width=0.75\textwidth,clip=]{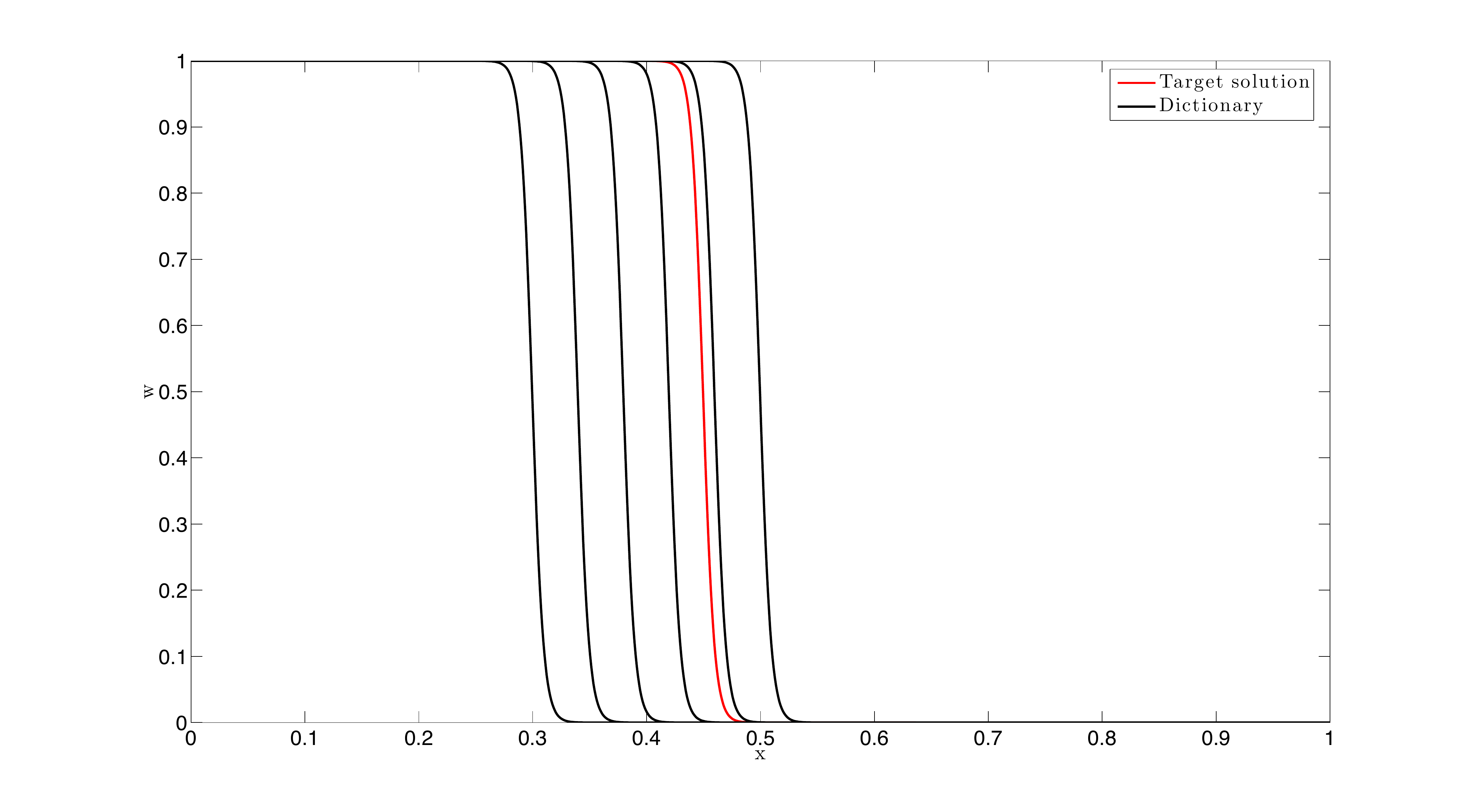}}
\end{center}
\caption{\label{fig:Advec1DDico}One-dimensional advection equation: dictionary and target solution}
\end{figure}

A dictionary of $r=6$ solutions is built for $\mu\in\{0.3,0.34,0.38,0.42,0.46,0.5\}$. The six solutions are depicted in Figure~\ref{fig:Advec1DDico}. Each solution has a high gradient  at a different location $x=\mu$. Five different model reduction methods are then compared, namely Galerkin projection, $L^1$-norm minimization by LP and IRLS, $L^2$-norm and Huber function minimizations.

\begin{figure}
\begin{center}
{\includegraphics[width=0.75\textwidth,clip=]{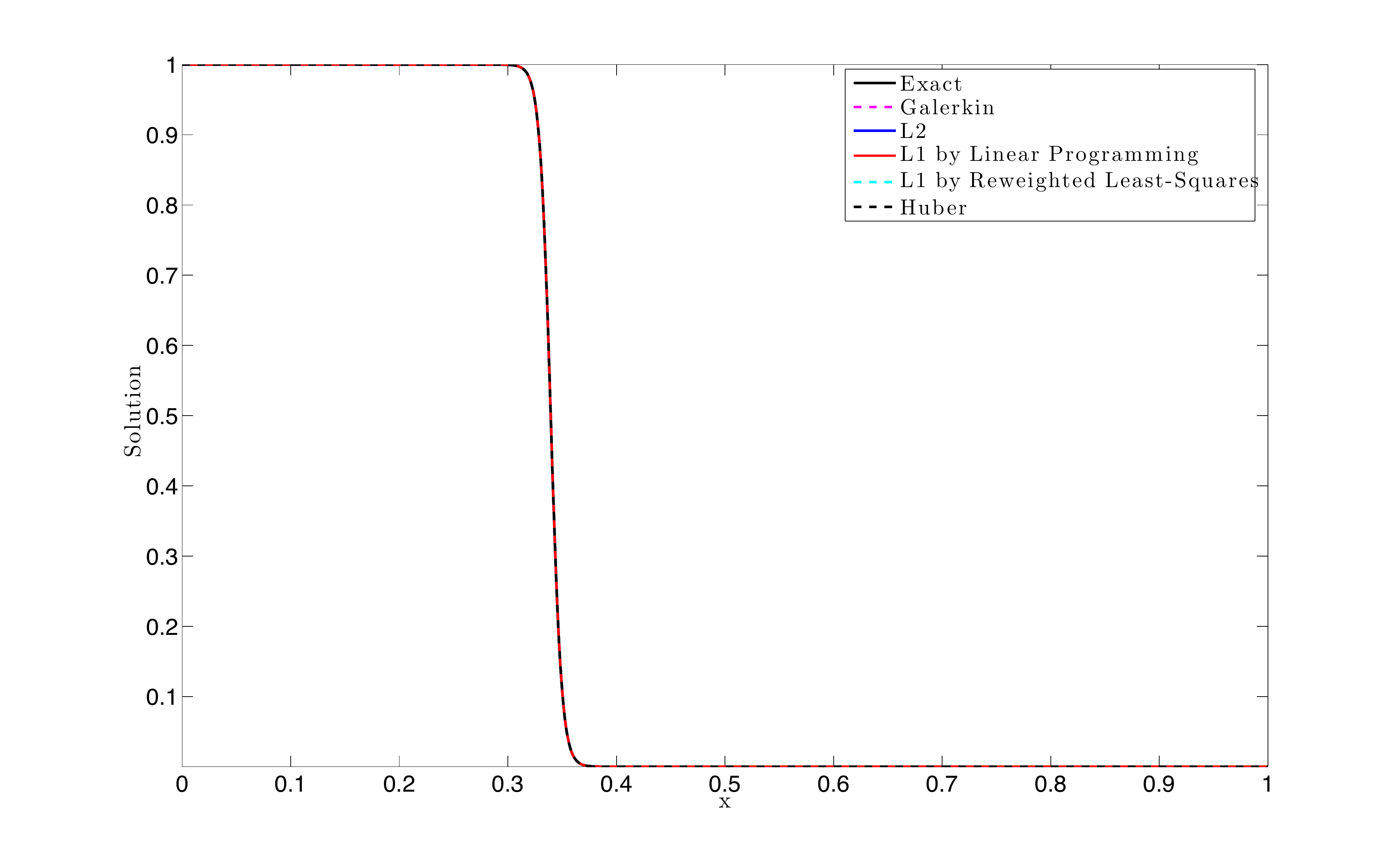}}
\end{center}
\caption{\label{fig:Advec1DDicoSolution}One-dimensional advection equation: solutions at one of the dictionary parameters $\mu^\star=0.34$} 
\end{figure}

\begin{figure}
\begin{center}
{\includegraphics[width=0.75\textwidth,clip=]{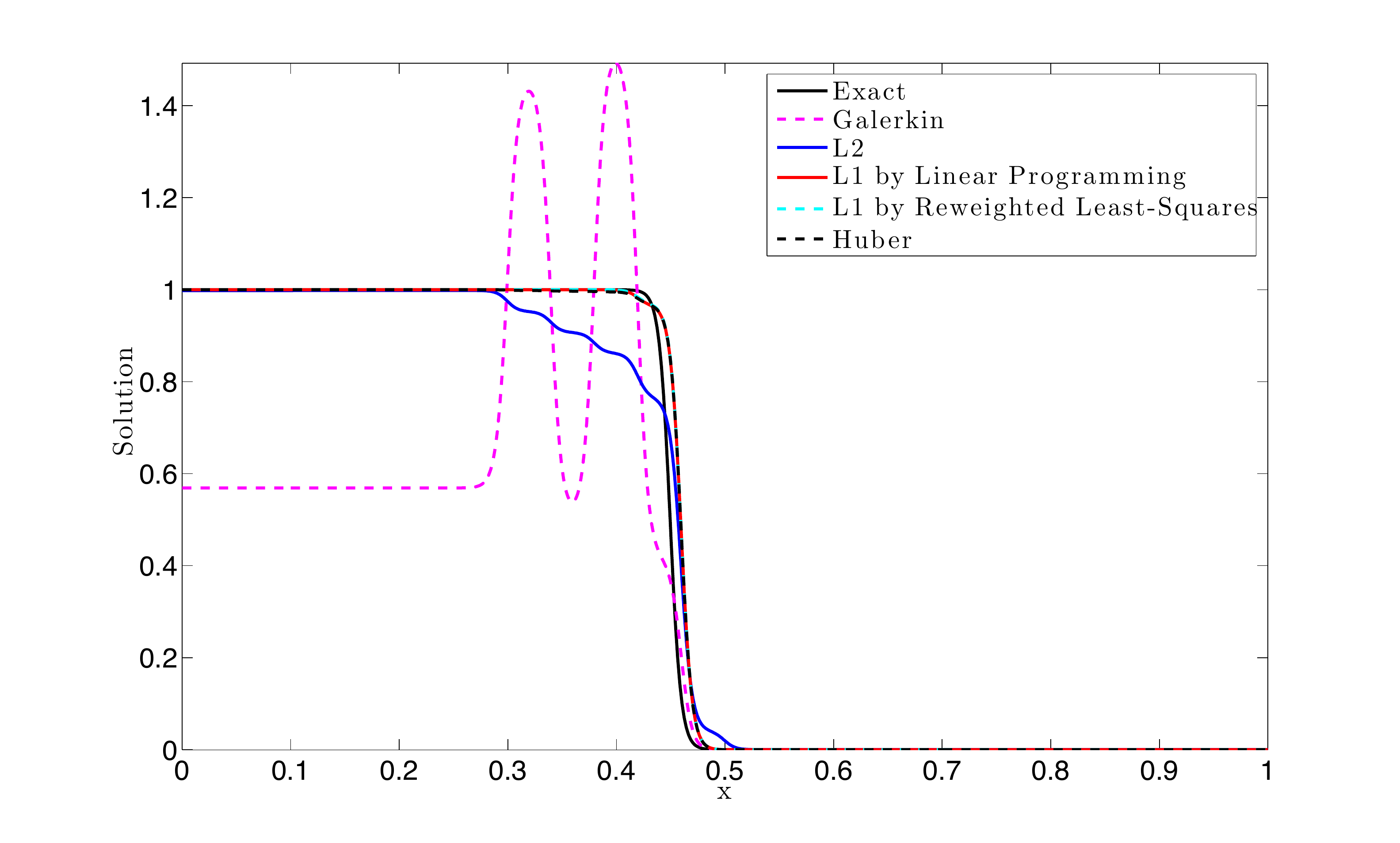}}
\end{center}
\caption{\label{fig:Advec1DSolution}One-dimensional advection equation: predicted solutions at target parameter $\mu^\star=0.45$}
\end{figure}

In a first experiment, a target parameter $\mu^\star=0.34$ belonging to the dictionary is selected. As reported in Figure~\ref{fig:Advec1DDicoSolution}, all approaches correctly predict the solution. Next, a target parameter $\mu^\star=0.45$ not belonging to the dictionary is selected. The model reduction results are reported in Figure~\ref{fig:Advec1DSolution}. One can observe that the $L^1$-norm and Huber function minimizations lead to predictions that are identical and are the most physical as they exhibit a single shock. On the other hand,  model reduction based on Galerkin projection predicts an unphysical solution with very large oscillations. $L^2$-norm minimization leads to a smooth solution that is not physical either. The reduced coordinates corresponding to each of the six model reduction approaches are reported in Table~\ref{tab:Advection1D}. One can observe that the solutions returned by the approaches based on  $L^1$-norm and Huber function minimizations are very similar and are sparse, unlike the solutions returned by Galerkin projection and $L^2$-norm minimization. Finally, the residuals returned by each minimization technique are depicted in Figure~\ref{fig:Advection1DResiduals}. One can observe that $L^2$-norm minimization results in a higher maximum local residual as well as non-zero residuals that have a much larger support in the computational domain.

 \begin{table}[htdp]
\begin{center}\begin{tabular}{|c|c|c|c|c|c|}
\hline
 & Galerkin  & $L^2$-norm  & $L^1$-norm  (LP) & $L^1$-norm (IRLS) & Huber function  \\
 \hline
 $\alpha_1$ &-0.896 & 0.046 & -3.6$\times 10^{-12}$ & -2.5$\times 10^{-11}$ & -1.4$\times 10^{-10}$ \\
 \hline
  $\alpha_2$  &0.962 &0.045 & 5.0$\times 10^{-8}$& 4.3$\times 10^{-8}$ & 1.4$\times 10^{-8}$\\
  \hline
    $\alpha_3$  &-1.028 &0.045 & -2.0$\times10^{-5}$ & -1.8$\times10^{-5}$ & -7.4$\times10^{-6}$\\
  \hline
    $\alpha_4$  &1.115 &0.097 & 0.033 & 0.031 & 0.019\\
  \hline
    $\alpha_5$  &0.417 &0.725 & 0.967 & 0.970 & 0.981\\
  \hline
    $\alpha_6$  &-0.001 &0.040 & -4.6$\times 10^{-4}$ & -5.0$\times 10^{-4}$ & -2.9$\times 10^{-4}$ \\
  \hline
 \end{tabular} 
 \caption{One-dimensional advection equation: reduced solutions}
 \label{tab:Advection1D}
\end{center}
\end{table}

\begin{figure}
\begin{center}
{\includegraphics[width=0.75\textwidth,clip=]{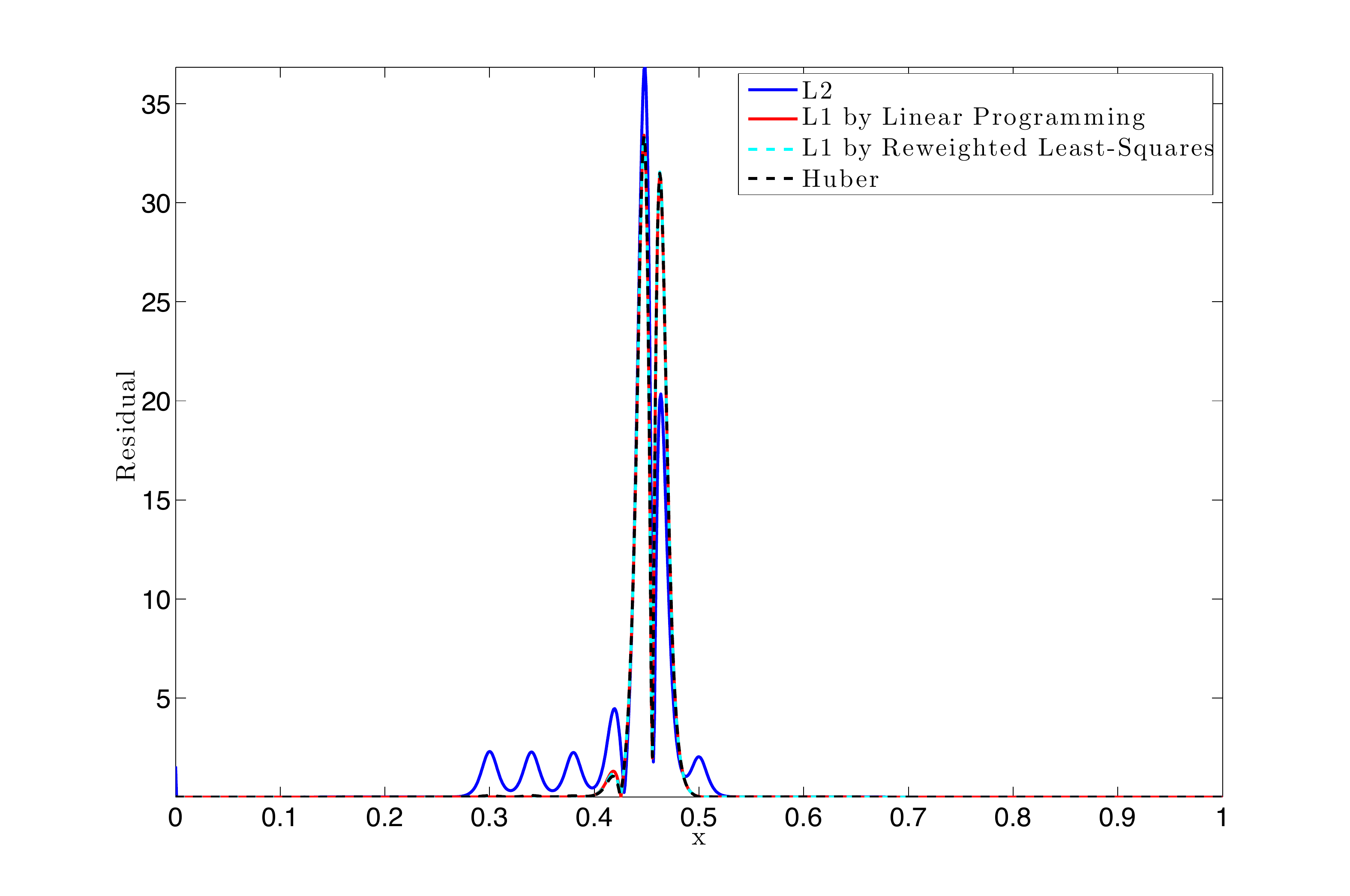}}
\end{center}
\caption{\label{fig:Advection1DResiduals}One-dimensional advection equation: residuals at target parameter $\mu^\star=0.45$}
\end{figure}
}


\subsubsection{{Two-dimensional advection-diffusion equation}}\label{ssec:advec2D}
\begin{figure}
\begin{center}
 \begin{subfigmatrix}{2}
  \subfigure[Galerkin projection]{ \includegraphics[width=5.in]{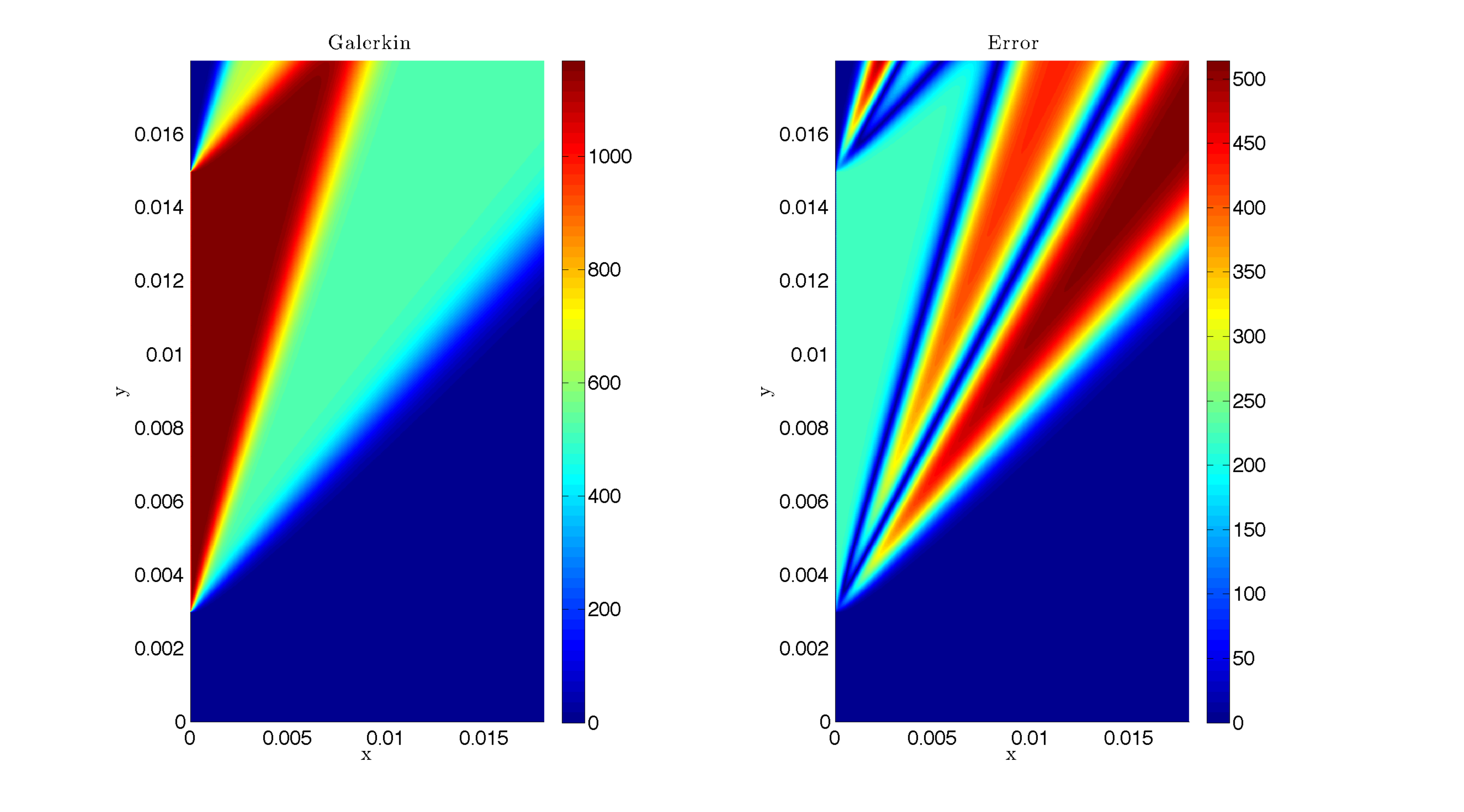}}
    \subfigure[$L^2$-norm minimization]{ \includegraphics[width=5.in]{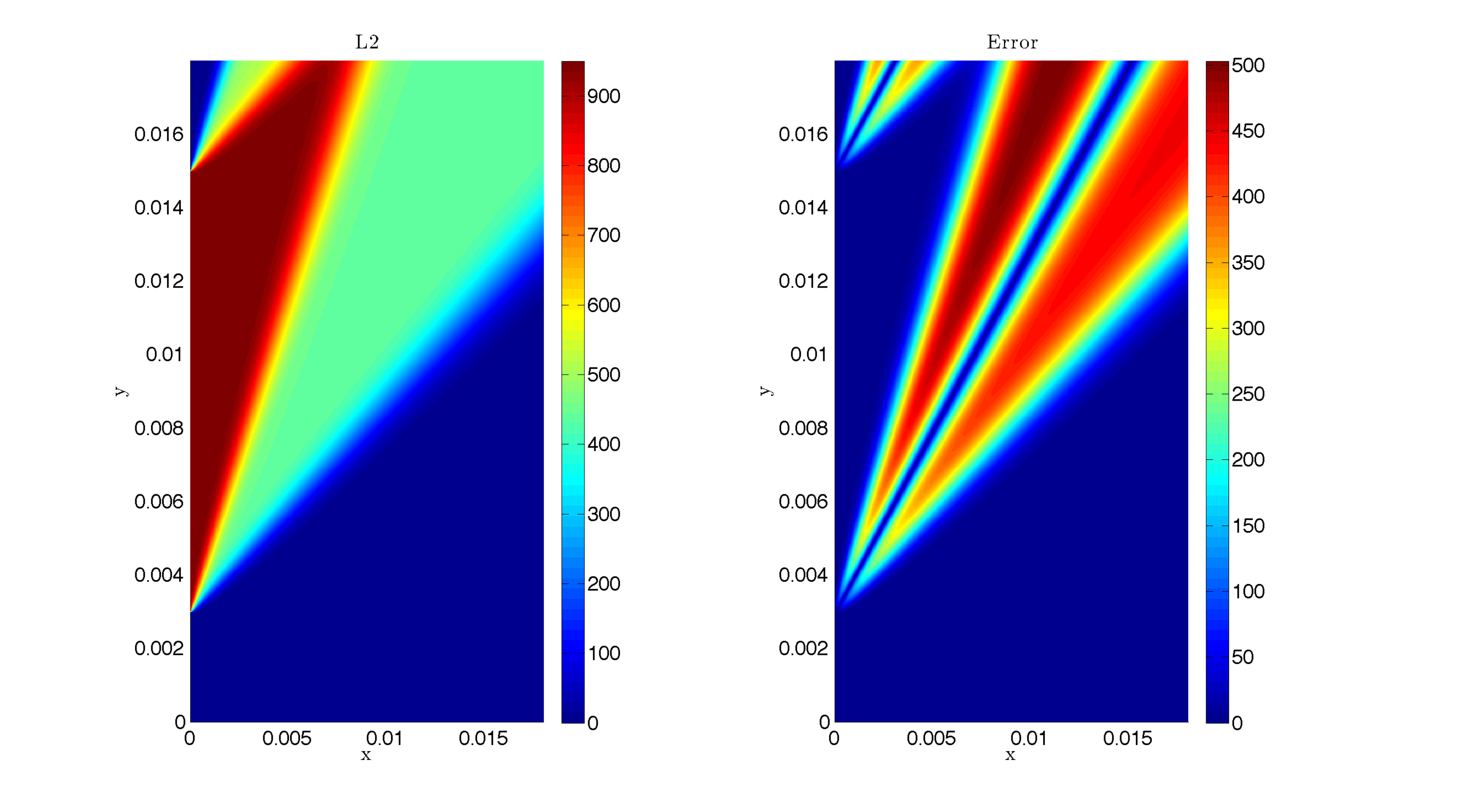}}
  \end{subfigmatrix}
 \caption{Two-dimensional advection-diffusion equation: predicted solutions and errors at target parameter $\mu^\star=\frac{\pi}{4}$}
 \label{fig:Advection2DSolution1}
 \end{center}
\end{figure}

\begin{figure}
\begin{center}
 \begin{subfigmatrix}{2}
  \subfigure[$L^1$-norm minimization]{ \includegraphics[width=5in]{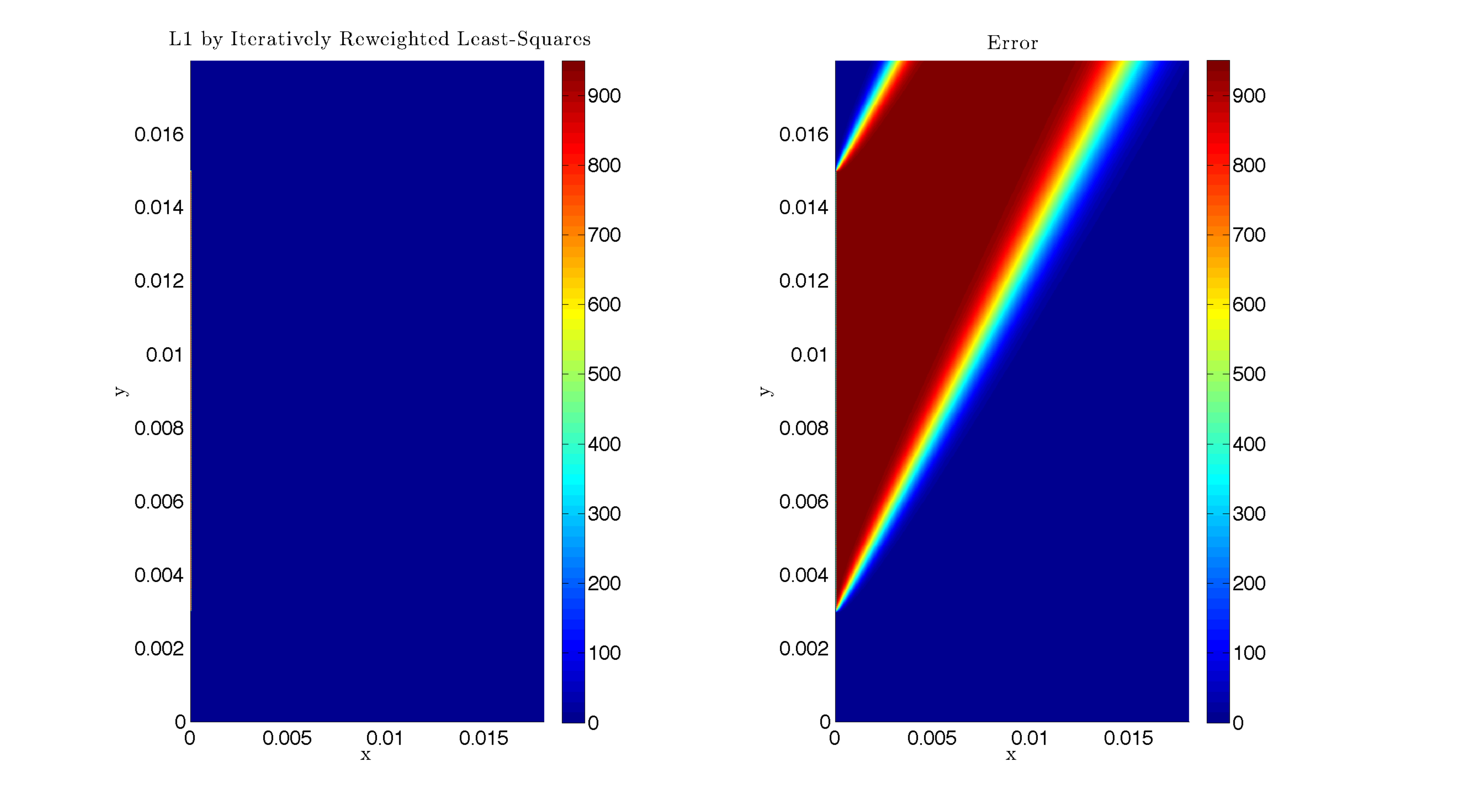}}
    \subfigure[Huber function minimization]{ \includegraphics[width=5.in]{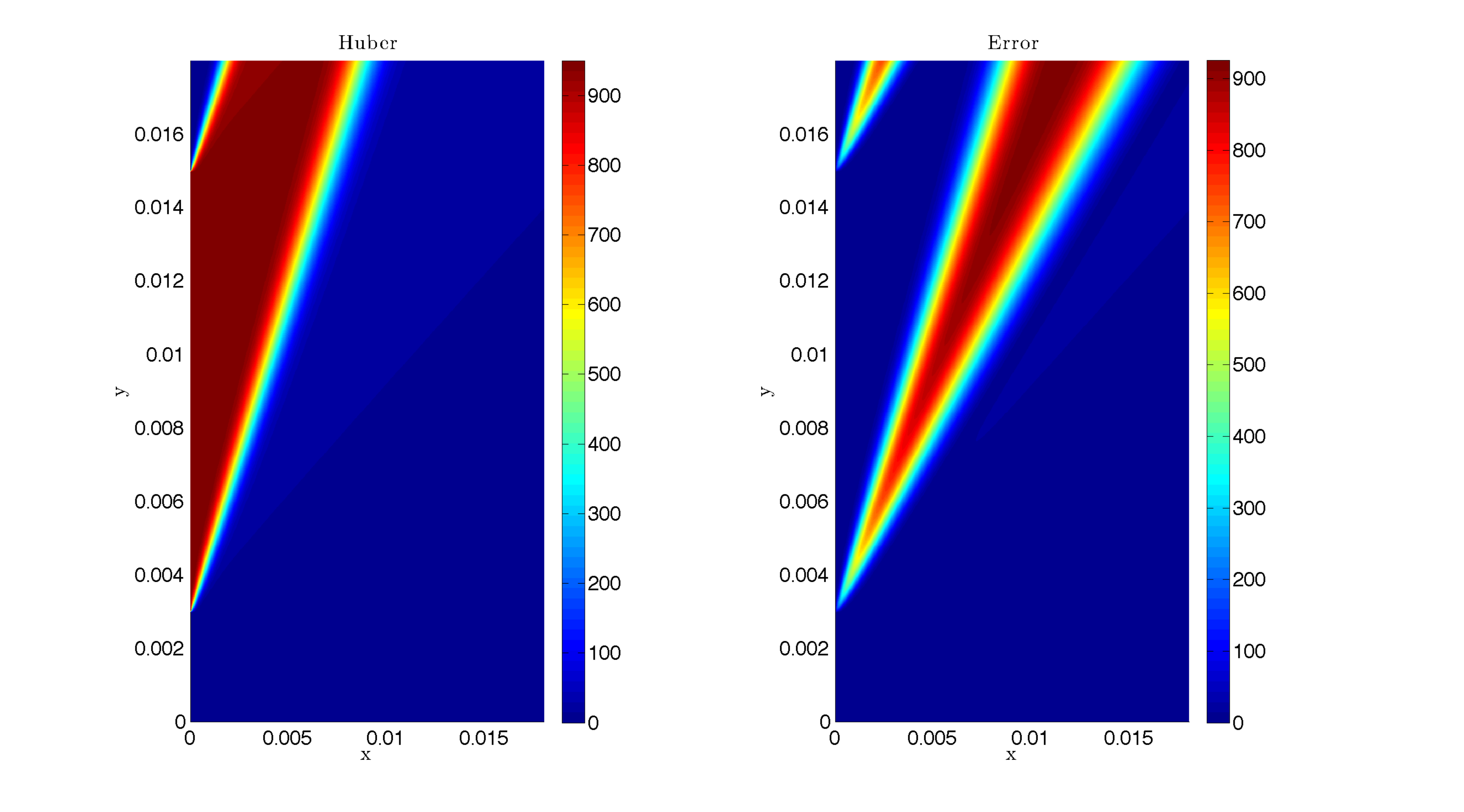}}
  \end{subfigmatrix}
 \caption{Two-dimensional advection-diffusion equation: predicted solutions and errors at target parameter $\mu^\star=\frac{\pi}{4}$ (continued)}
 \label{fig:Advection2DSolution2}
 \end{center}
\end{figure}

{The two-dimensional advection-diffusion equation is then considered
\begin{equation}
\lambda(\mu) \cdot\nabla  u(x,y) -\kappa  \Delta u(x,y)  = 0,~(x,y) \in \Omega = [0,0.018] \times  [0,0.018]
\end{equation}
with incoming Dirichlet boundary conditions and outgoing Neumann boundary conditions. The problem is parameterized by the angle  of the advection flow with respect to the $x$ axis: $\lambda(\mu) = (\|\lambda\|_2\cos(\mu),\|\lambda\|_2\sin(\mu))$. This problem is dominated by  advection since $\|\lambda\|_2=0.5$ and $\kappa=2\times 10^{-7}$.  

The problem is discretized by finite differences using a uniform mesh with $304$ points in each direction, resulting in $N=88464$ degrees of freedom. For this large scale problem, solving the $L^1$-norm minimization problem by LP is not tractable and as such only the IRLS method is used in the $L^1$-norm case.

A dictionary $\mathcal{D}$ of two solutions is constructed for $\mu\in\{\frac{\pi}{6},\frac{\pi}{3}\}$ and a target parameter $\mu^\star=\frac{\pi}{4}$ considered. The respective solutions are depicted in Figure~\ref{fig:Advection2DDico}.

\begin{figure}
\begin{center}
 \begin{subfigmatrix}{2}
  \subfigure[Dictionary member \#1, $\mu=\frac{\pi}{3}$]{ \includegraphics[width=3.in]{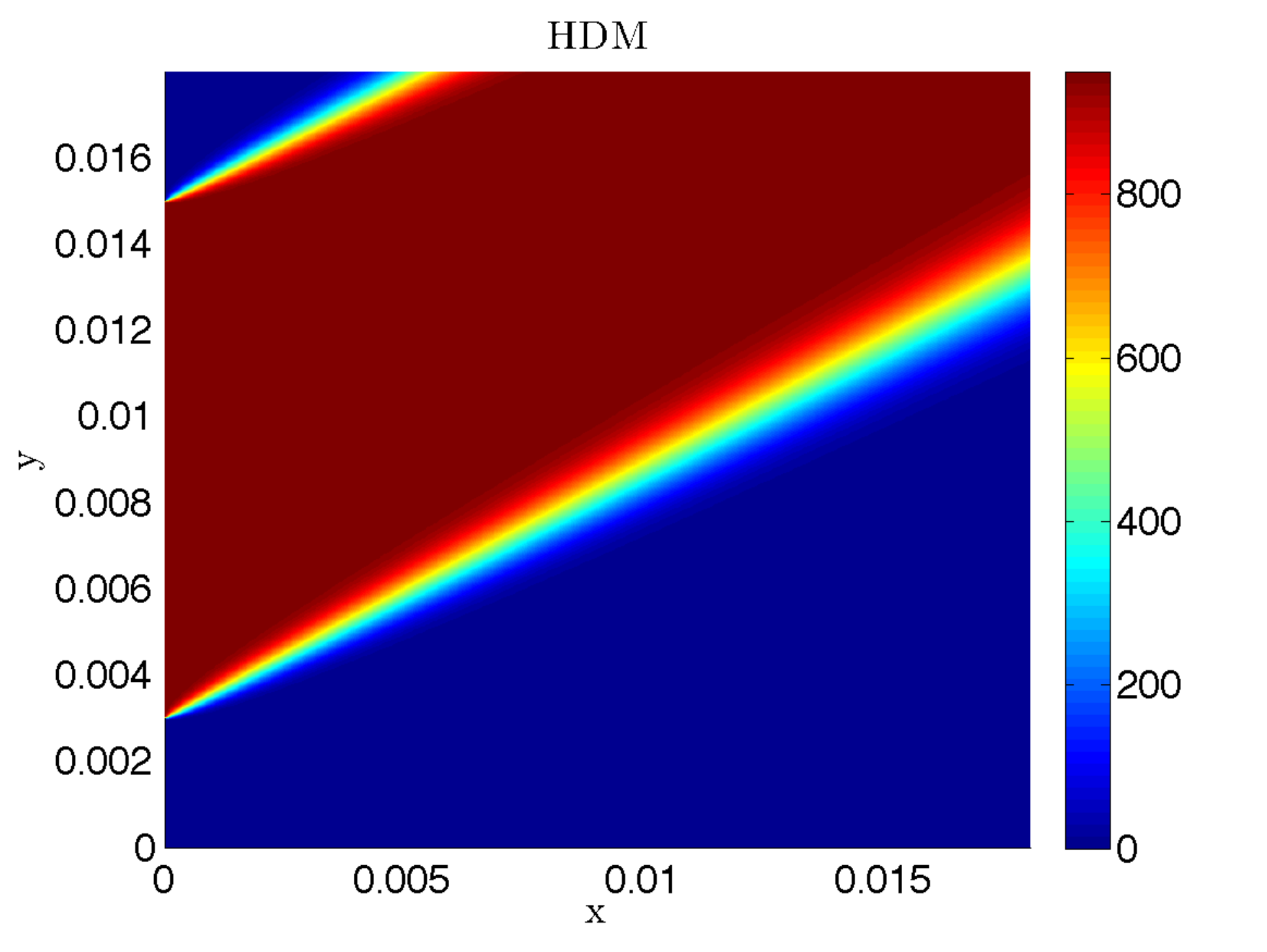}}
    \subfigure[Dictionary member \#2, $\mu=\frac{\pi}{6}$]{ \includegraphics[width=3.in]{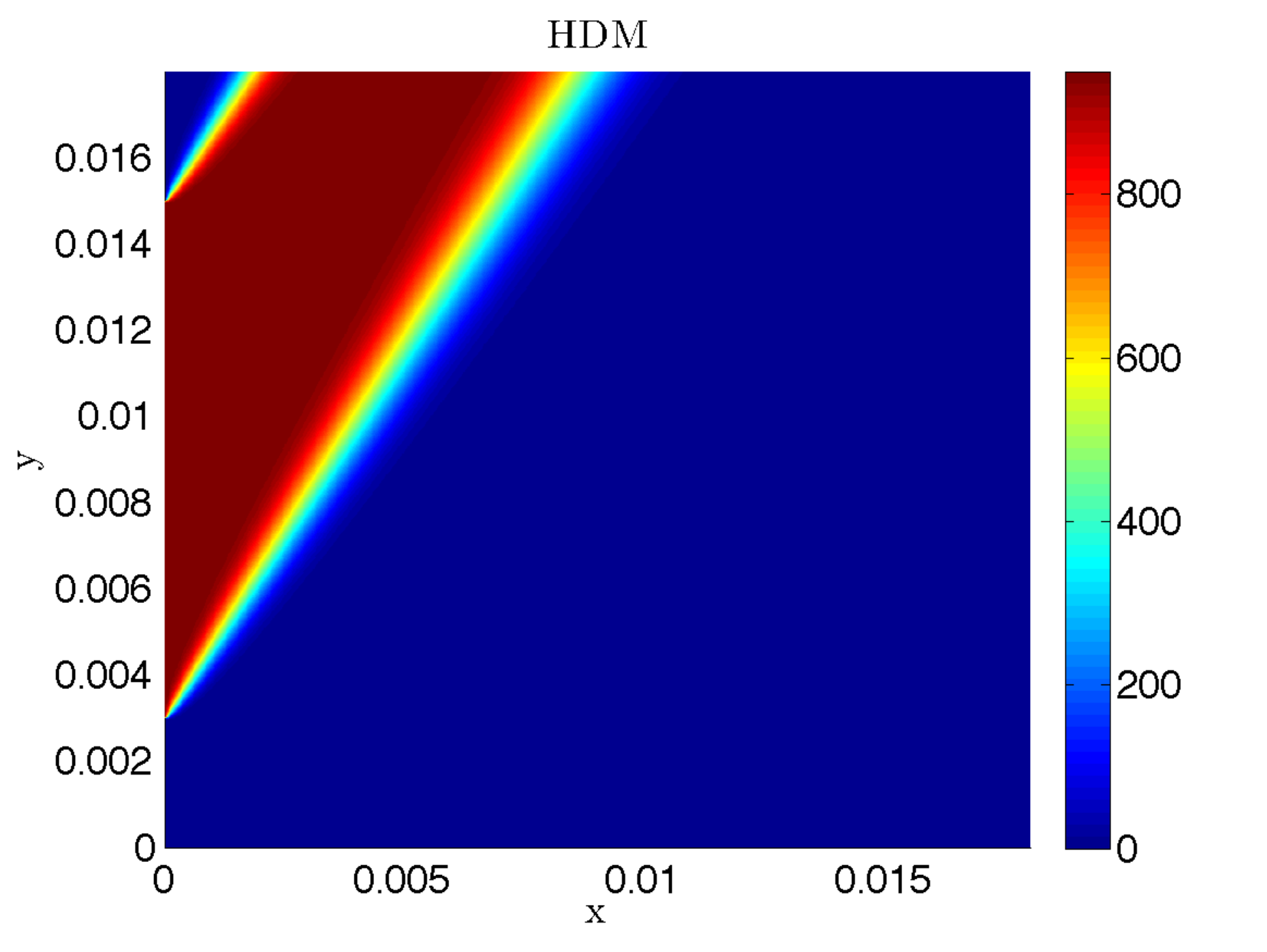}}
  \subfigure[Target solution $\mu^\star=\frac{\pi}{4}$]{ \includegraphics[width=3.in]{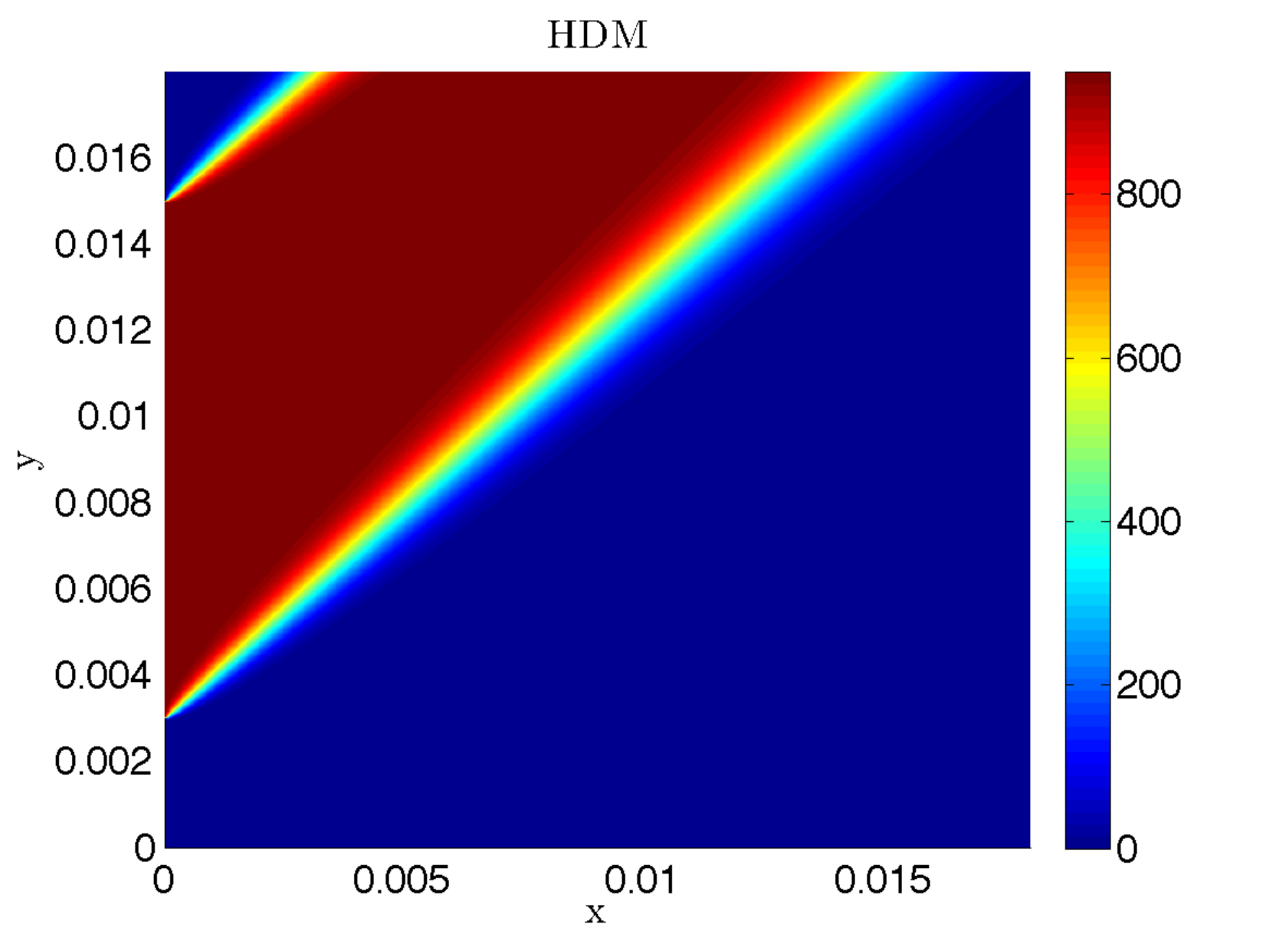}}
  \end{subfigmatrix}
 \caption{Two-dimensional advection-diffusion equation: dictionary and target solutions}
 \label{fig:Advection2DDico}
 \end{center}
\end{figure}

The following four model reduction methods are then applied: Galerkin projection, $L^2$-norm minimization, $L^1$-norm minimization by IRLS and Huber function minimization. The corresponding reduced solutions are reported in Table~\ref{tab:Advection2D} and the solutions and errors in Figures~\ref{fig:Advection2DSolution1}--\ref{fig:Advection2DSolution2}. For this problem, the $L^1$-norm method by IRLS failed to converge and the returned solution is zero. However, the Huber function minimization approach was much more robust and returned a physical solution with sharp gradient. Hence, this example illustrates the advantage of using the Huber function versus pure $L^1$-norm minimization.  Galerkin projection and $L^2$-norm minimization returned  very similar but much less physical solutions with gradients that are much less sharp.

 \begin{table}[htdp]
\begin{center}\begin{tabular}{|c|c|c|c|c|}
\hline
 & Galerkin  & $L^2$-norm  &  $L^1$-norm (IRLS) & Huber function  \\
 \hline
 $\alpha_1$ & 0.543 &  0.467 & 1.9$\times 10^{-11}$ & 0.021\\
 \hline
  $\alpha_2$  & 0.688 &  0.529& 5.2$\times 10^{-10}$ &  0.979\\
  \hline
 \end{tabular} 
 \caption{Two-dimensional advection-diffusion equation: reduced solutions}
 \label{tab:Advection2D}
\end{center}
\end{table}

}

\subsubsection{{Steady Burgers' equation}}

\begin{figure}
\begin{center}
{\includegraphics[width=0.75\textwidth,clip=]{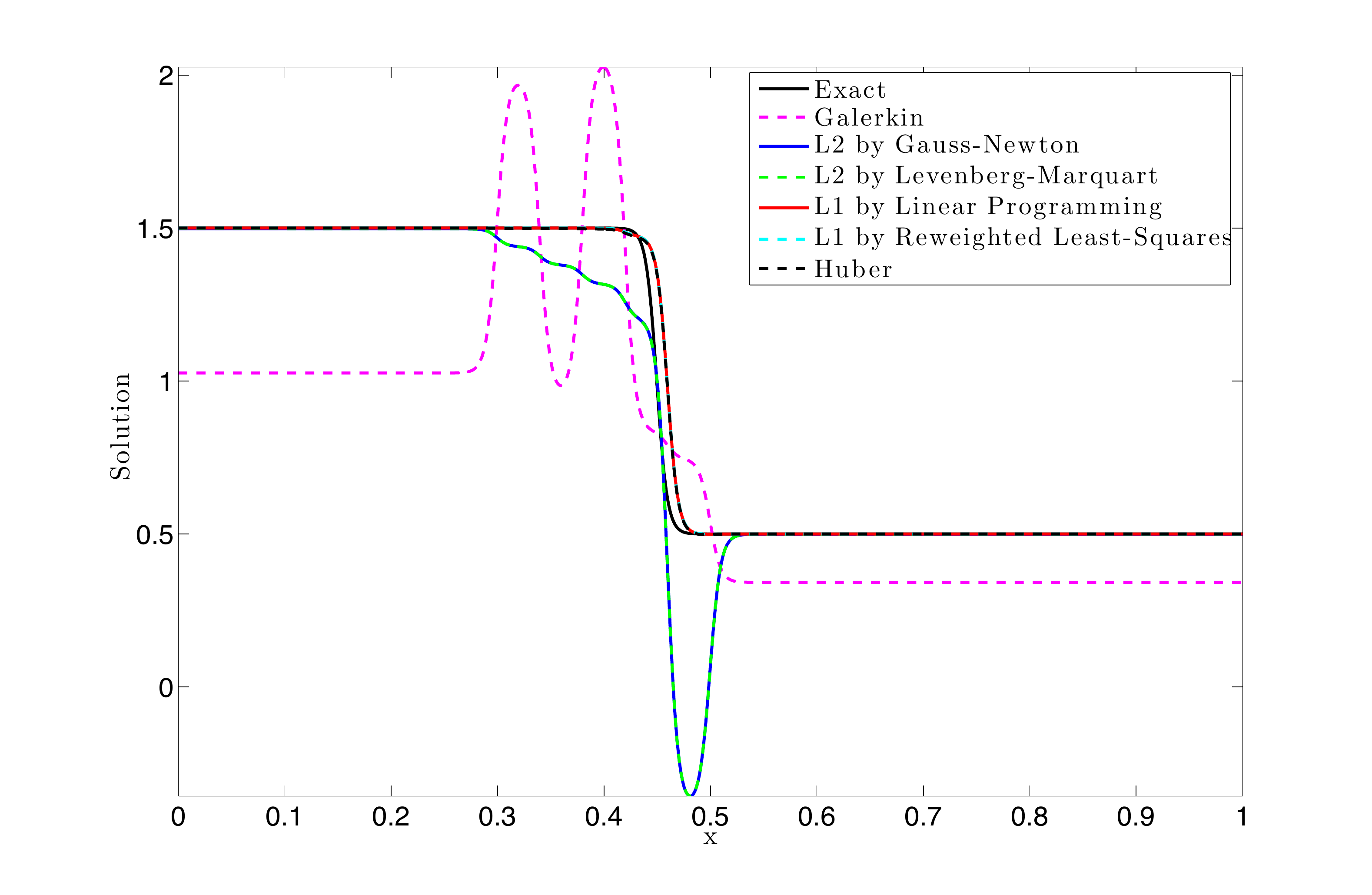}}
\end{center}
\caption{\label{fig:Burgers1DSolution}One-dimensional Burgers' equation: predicted solutions at target parameter $\mu^\star=0.45$}
\end{figure}

\begin{figure}
\begin{center}
{\includegraphics[width=0.75\textwidth,clip=]{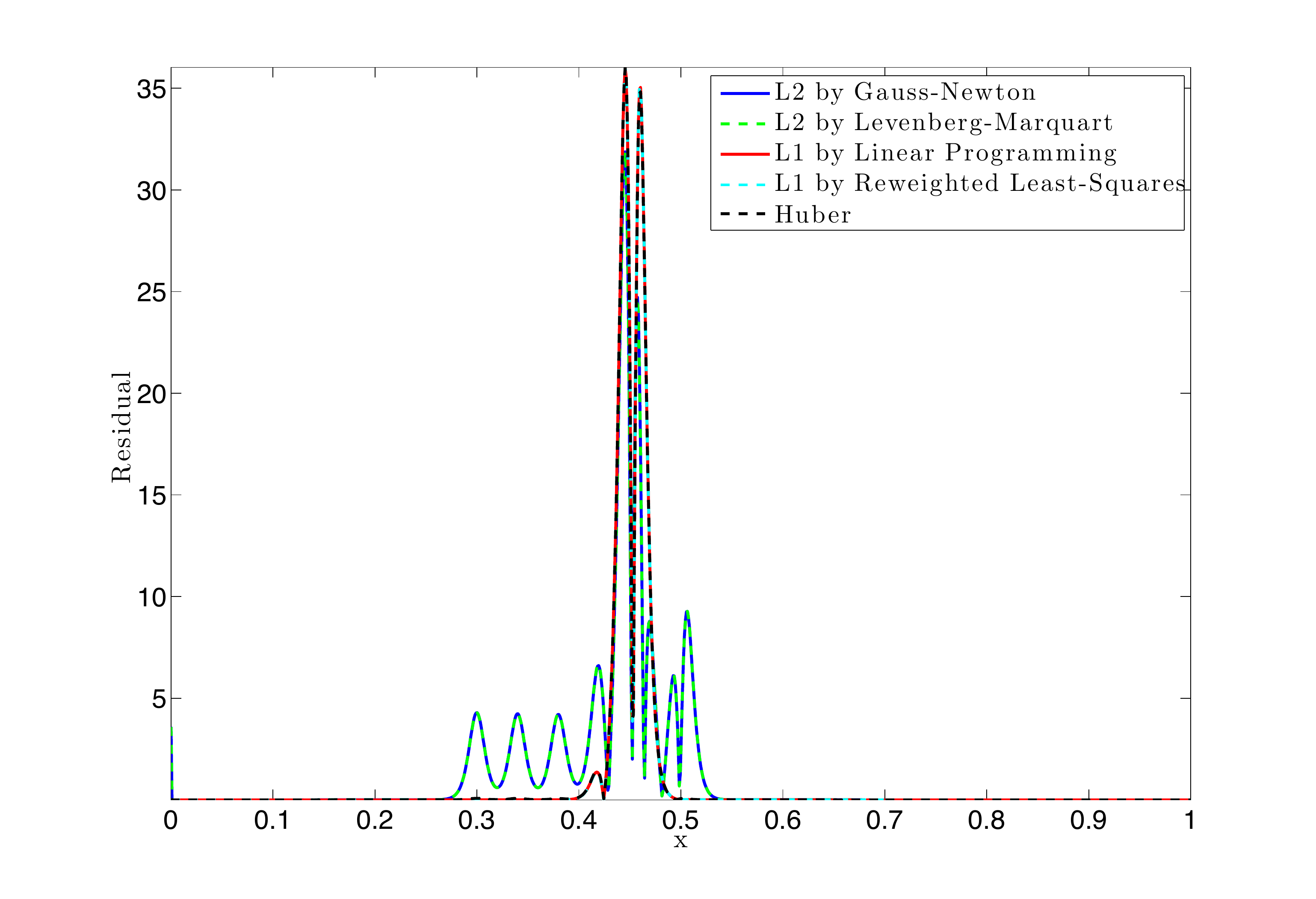}}
\end{center}
\caption{\label{fig:Burgers1DResiduals}One-dimensional Burgers' equation: residuals at target parameter $\mu^\star=0.45$}
\end{figure}

{ The following one-dimensional steady Burgers' equation is considered:
\begin{equation}
\frac{1}{2}\dpar{(u^2)}{x}(x)
= f(x;\mu),~x\in\Omega = [0,~1], 
\end{equation}
where $$f(x;\mu) = \frac{-2k\exp(-2k(x-\mu))(1+3\exp(-2k(x-\mu))}{(1+\exp(-2k(x-\mu)))^3}$$ with $k=100$. The solution exhibits a strong gradient  at location $x=\mu$. A Dirichlet boundary condition $u(0) = 1.5$ is applied at $x=0$. This PDE is discretized by finite differences using a uniform mesh, resulting in a HDM of dimension $N=10^3$.

A dictionary of $r=6$ solutions is built for $\mu\in\{0.3,0.34,0.38,0.42,0.46,0.5\}$ by solving each steady state problem by Newton-Raphton's method. The six solutions are depicted in Figure~\ref{fig:Burgers1DDico}. Each solution has a ``shock" at a different location. For this case, six different model reduction methods are compared, namely Galerkin projection, $L^2$-norm minimization by Gauss-Newton and Levenberg-Marquardt, $L^1$-norm minimization by LP and IRLS and Huber function minimization.

\begin{figure}
\begin{center}
{\includegraphics[width=0.75\textwidth,clip=]{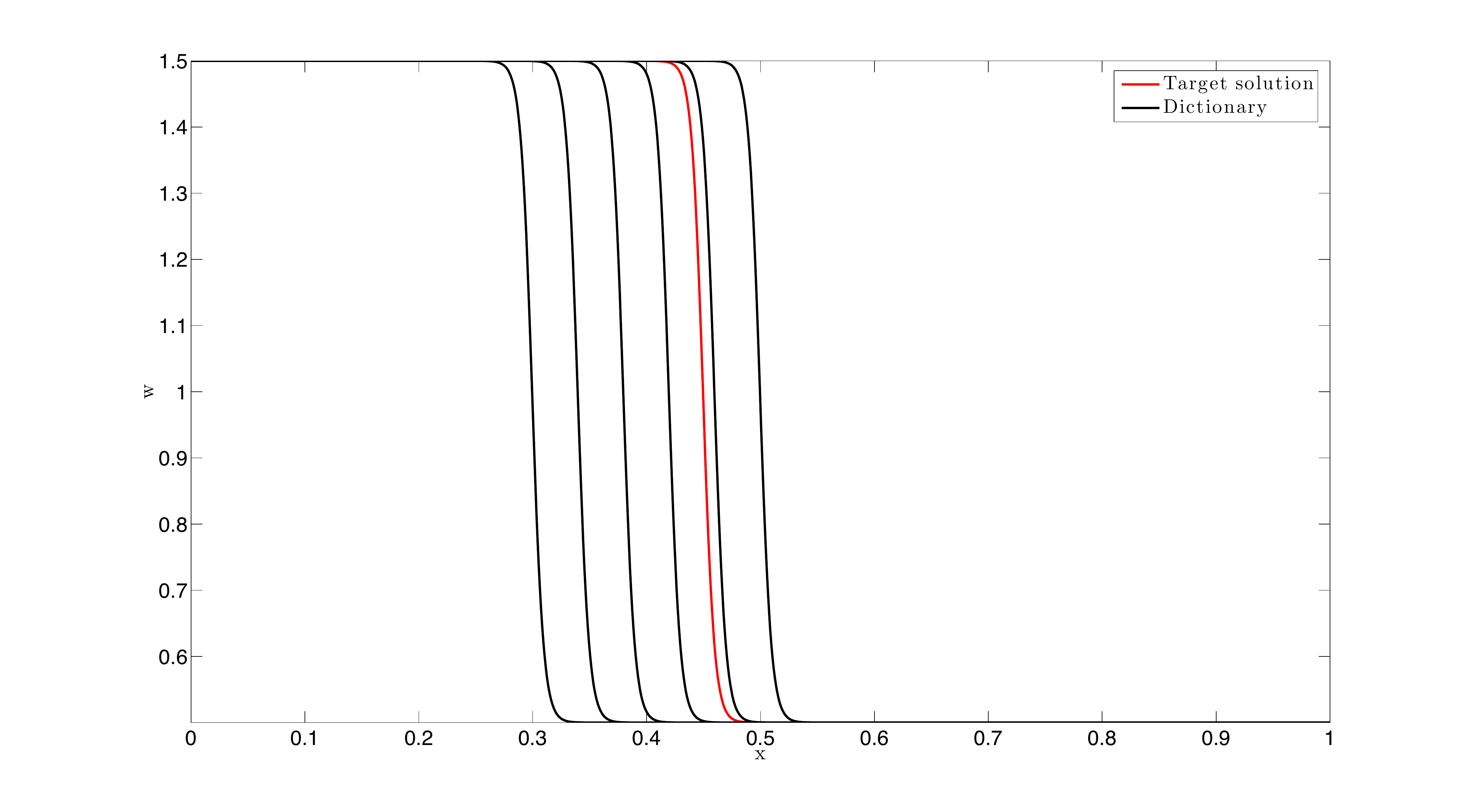}}
\end{center}
\caption{\label{fig:Burgers1DDico}One-dimensional Burgers' equation: dictionary and target solution}
\end{figure}

A target parameter $\mu^\star=0.45$ not belonging to the dictionary is selected. The model reduction results are reported in Figure~\ref{fig:Burgers1DSolution}. One can observe that the $L^1$-norm and Huber function minimization results are identical and are the most physical as they exhibit a single ``shock". On the other hand,  Galerkin projection results in an unphysical solution with very large oscillations and inaccurate constant solutions before and after the ``shock".  The two $L^2$-norm-based approaches result in a very smooth solution before the shock and an undershoot after the shock that are not physical either.  Finally, the residuals returned by each minimization technique are shown in Figure~\ref{fig:Burgers1DResiduals}. One can observe that $L^2$-norm minimization approaches result in a higher maximum local residual as well as non-zero residuals that have a larger support in the computational domain.

}

\subsection{Model reduction of unsteady problems}
\subsubsection{Unsteady Burgers' equation}\label{Unsteady Burgers}
\begin{figure}[ht]
\subfigure[Solutions]{\includegraphics[width=0.45\textwidth]{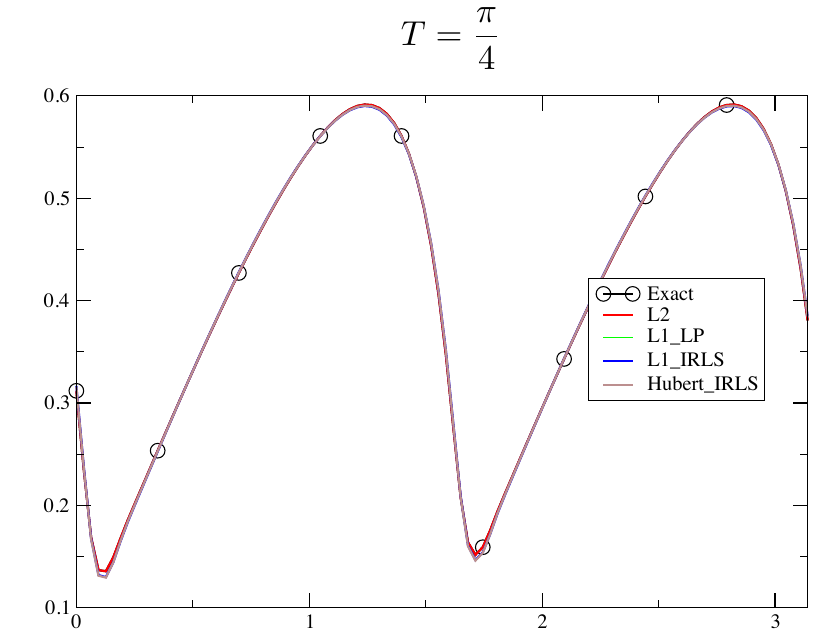}}
\subfigure[Zoom near a maximum]{\includegraphics[width=0.45\textwidth]{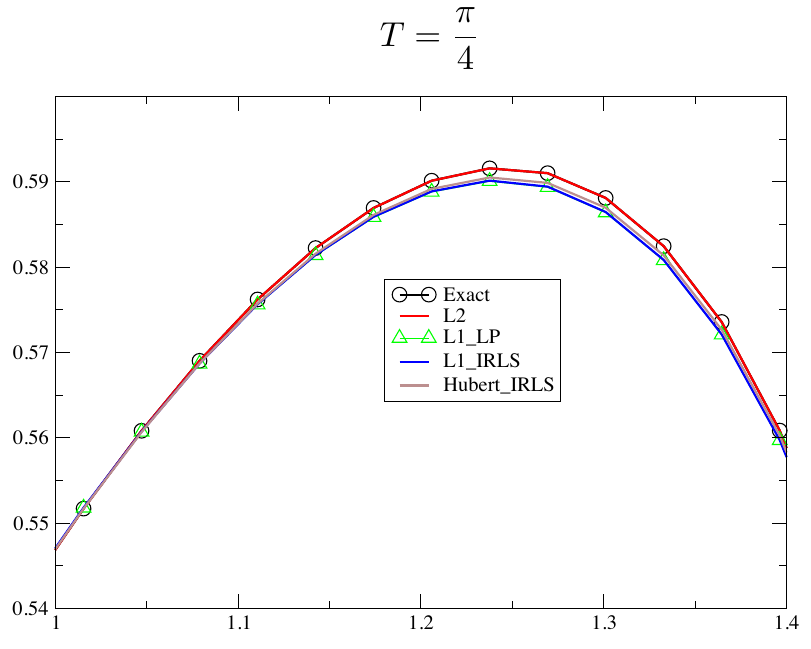}}
\centering{\subfigure[Zoom near a minimum]{\includegraphics[width=0.45\textwidth]{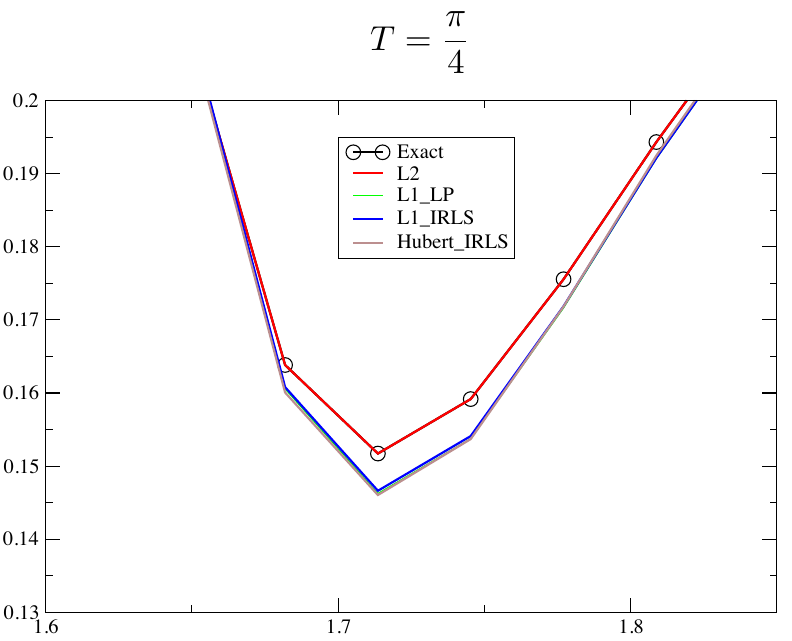}}}
\caption{\label{burger:pis4} Unsteady Burgers' equation: predicted solutions at target parameter $\mu^\star=0.5$ at $t=\frac{\pi}{4}$ }
\end{figure}

\begin{figure}[ht]
\subfigure[]{\includegraphics[width=0.45\textwidth]{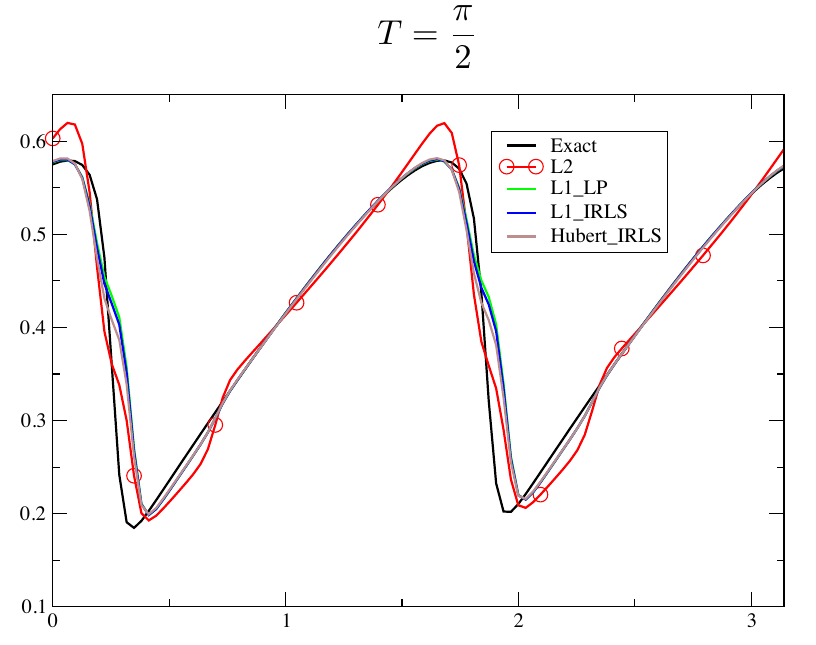}}
\subfigure[]{\includegraphics[width=0.45\textwidth]{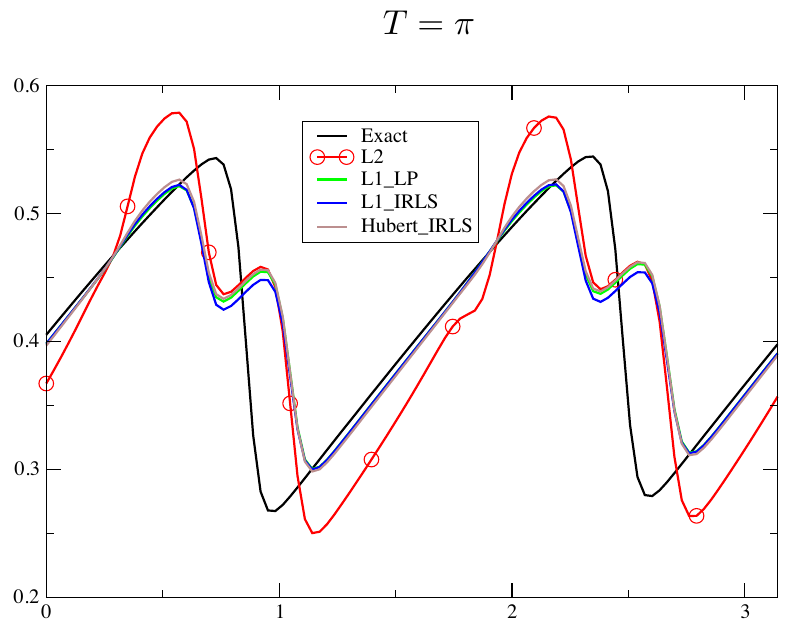}}
\caption{\label{burger:apres_choc} Unsteady Burgers' equation: predicted solutions at target parameter $\mu^\star=0.5$ at $t=\frac{\pi}{2}$ (left) and $t=\pi$ (right)}
\end{figure}

We consider here the system~\eqref{eq:2} in $\Omega=[0,2\pi]$ with periodic boundary conditions and initial conditions parameterized  by
$$
u_0(x;\mu)=\mu \; \big | \sin(2\; x)\big |+0.1,
$$
where $\mu \in [0,1]$. In this setting, the solution develops a shock that travels with the velocity $\sigma_\mu=0.55 \mu$.
A dictionary $\mathcal{D}$ is constructed by sampling the parameters $\{0,0.2,0.4,0.6,1.0\}$ ($r=5$) and the solution sought for the predictive case $\mu^\star=0.5$. A shock appears at $t=1$. We display the solutions obtained by $L^2$-norm, $L^1$-norm by LP and IRLS and the Hubert-IRLS minimization procedures for $t=\tfrac{\pi}{4}<1$, $t=\tfrac{\pi}{2}$ and $t=\pi$ in Figures \ref{burger:pis4} and \ref{burger:apres_choc}.
Before the shock appears, there is almost no difference between the four solutions. The approach based on minimizing the $L^2$-norm is even slightly better, as it can be observed  from the two zooms in Figure~\ref{burger:pis4}. However, the situation after the shock appears is very different, as observed in Figure \ref{burger:apres_choc}: the $L^2$-norm solution is clearly the worst one with large oscillations. The $L^1$-norm-type solutions are all close to each other and the shock is rather well reproduced with, however, an artifact that develops for longer  times, as seen at $t=\pi$. Nevertheless, the $L^1$-norm-type solutions are within the bounds of the ``exact" solution, and no large oscillation develops.

In a second set of numerical experiments, we consider the influence of the sampling parameter set included in the dictionary $\mathcal{D}$. We consider two dictionaries
 $\mathcal{D}_1=\{0.4,0.45, 0.55$ $ , 0.6\}$ and  $\mathcal{D}_0=\{0,0.2,0.4,0.45, 0.55 , 0.6,1.0\}$, for the same target value of $\mu^\star=0.5$. These choices amounts to selecting samples close to the target value $0.5$ while varying elements of the dictionary that are not close to $0.5$.
\begin{figure}[ht]
\subfigure[$\mathcal{D}_1$]{\includegraphics[width=0.45\textwidth]{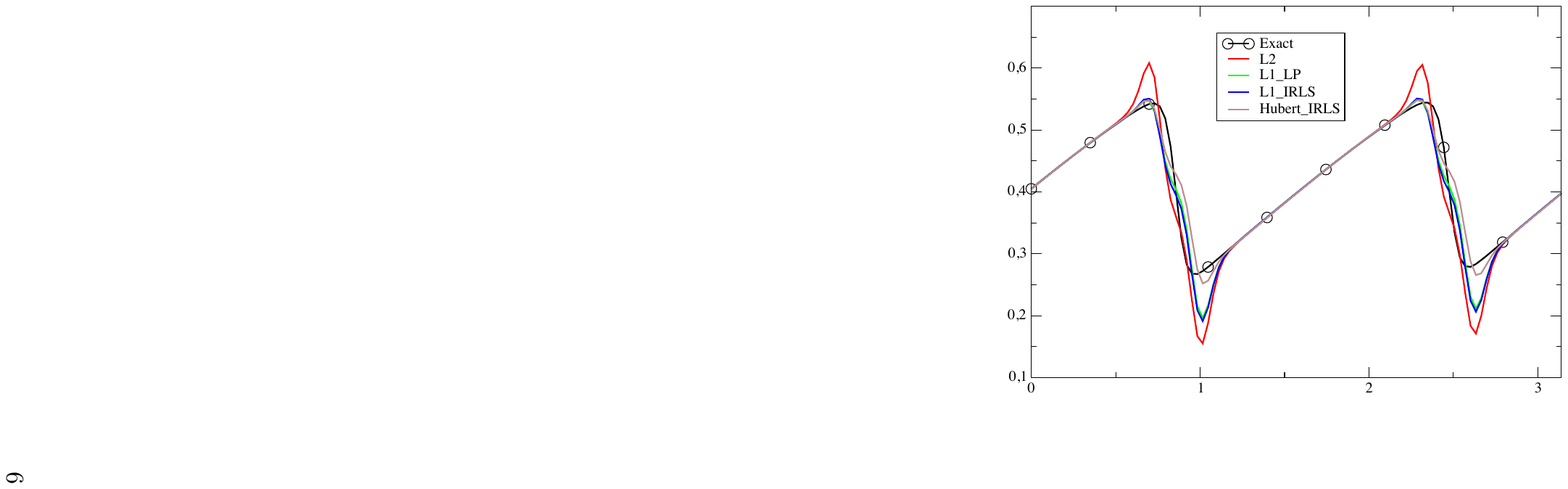}}
\subfigure[$\mathcal{D}_0$]{\includegraphics[width=0.45\textwidth]{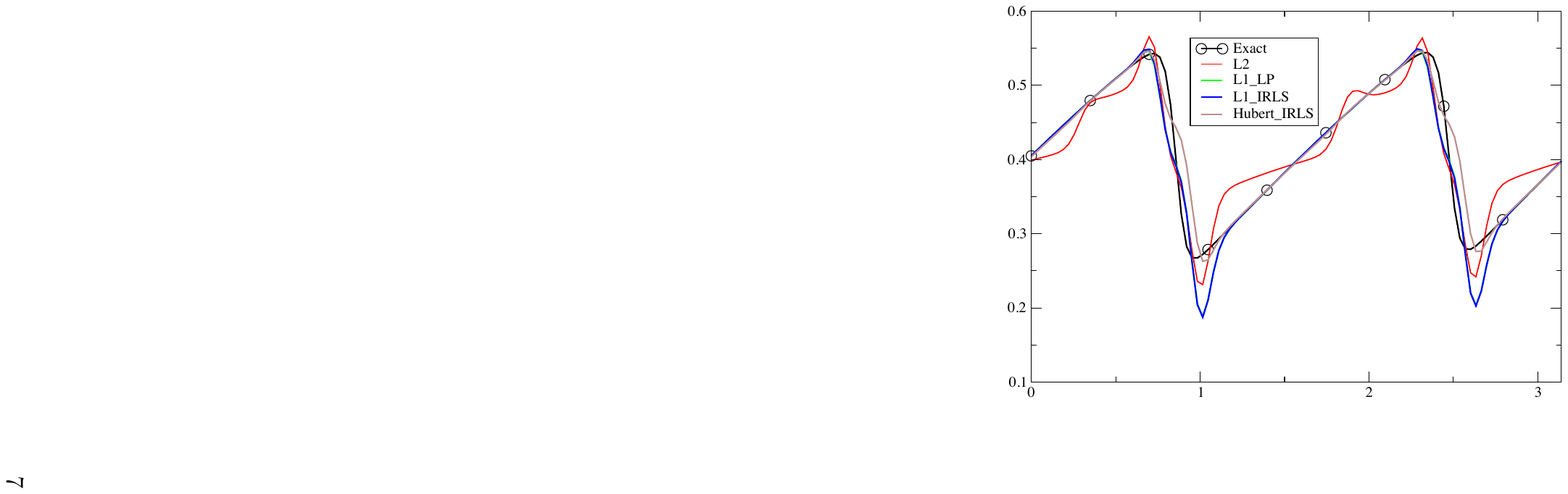}}
\caption{\label{influence} Unsteady Burgers' equation: predicted solutions at target parameter $\mu^\star=0.5$ at $t=\pi$ for two dictionaries associated with two samples of the parameter domain $\mathcal{P}$}
\end{figure}
We see that refining the dictionary has a positive influence as the target solution is much closer to the dictionary elements. This is confirmed by additional experiments where the samples of  $\mu$ used to generate the dictionary were more numerous and closer to $0.5$ (not reported here). However, keeping values of $\mu$ that are far from $0.5$ has an effect,  clearly negative for the $L^2$-norm optimization. The $L^1$-norm-type solutions are however  unaffected by the presence of these ``outliers" in the dictionary, similarly as in the simple regression case reported in Section~\ref{ssec:regression}. Overall, the most accurate predictions are  based on the minimization of  the Hubert function by IRLS.

%

\subsubsection{Euler equations}
The one-dimensional compressible Euler equations are considered on $\Omega = [0,1]$
\begin{subequations}\label{euler}
\begin{equation}\label{euler:1}
\dpar{}{t}\begin{pmatrix}\rho \\ \rho u \\E\end{pmatrix}+\dpar{}{x}\begin{pmatrix}
\rho u \\ \rho u^2+p \\ u(E+p)\end{pmatrix}=0,
\end{equation}
for which $U=(\rho,\rho u,E)^T$ and the pressure is given by
\begin{equation}
\label{euler:2}
p=(\gamma-1)\left( E-\frac{1}{2}\rho u^2\right)
\end{equation}
with $\gamma=1.4$.

This problem is parameterized by the initial conditions $U_0(x;\mu)$. To define the parameterized initial conditions of the problem, the Lax and Sod cases are first introduced as follows. 

The state  $U_{\text{Sod}}(x)$ is defined by the primal physical quantities:
\begin{equation}
\label{euler:4}
V_{\text{Sod}}(x) =\left\{\begin{array}{ll}
\rho= 1 \text{ if }x\leq 0.5, & 0.125 \text{ otherwise,}\\
u=0.0 \\
p=1.0 \text{ if }x\leq 0.5, & 0.1 \text{ otherwise,}
\end{array}\right .
\end{equation}
and $U_{\text{Lax}}(x)$ defined by
\begin{equation}
\label{euler:5}
V_{\text{Lax}}(x)=\left\{\begin{array}{ll}
\rho= 0.445 \text{ if }x\leq 0.5, & 0.5 \text{ otherwise,}\\
u=0.698 \text{ if }x\leq 0.5, & 0.0 \text{ otherwise,}\\
p=3.528 \text{ if }x\leq 0.5, & 0.571 \text{ otherwise.}
\end{array}\right .
\end{equation}

The Sod condition presents a fan, followed by a contact and a shock. For the density and the pressure, the solution behaves monotonically, and the contact is moderate. The Lax solution has a very different behavior and the contact is much stronger. This is depicted in Figure \ref{figure:sod-lax} where the two solutions are shown for $t=0.16$.
\begin{figure}[h]
\begin{tabular}{cc}
\includegraphics[width=0.45\textwidth]{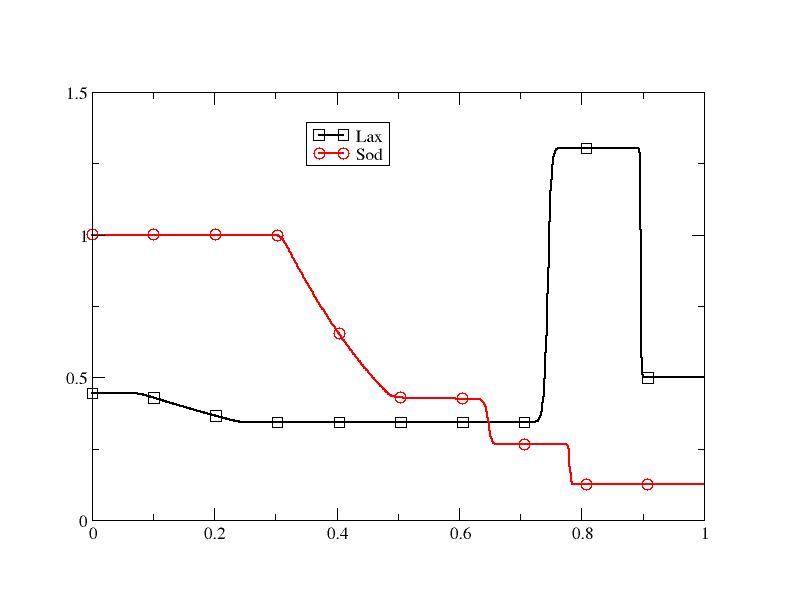}&  \includegraphics[width=0.45\textwidth]{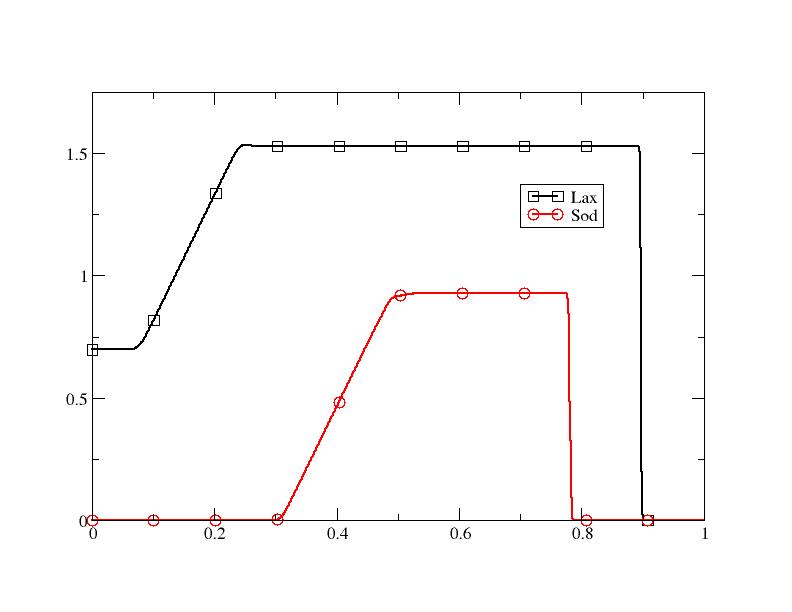}\\
Density& Velocity\\
\multicolumn{2}{c}{\includegraphics[width=0.45\textwidth]{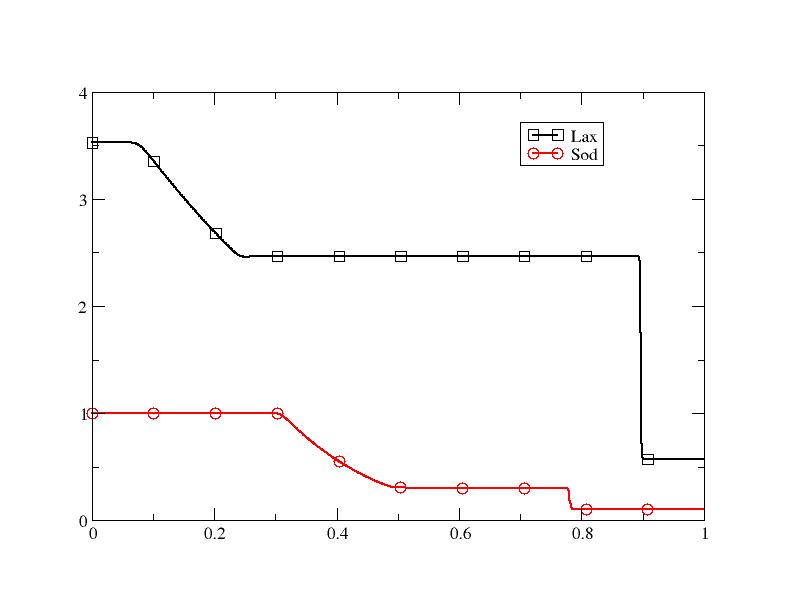}}\\
\multicolumn{2}{c}{Pressure}
\end{tabular}
\caption{\label{figure:sod-lax} One-dimensional Euler equations: density, velocity and pressure for the Lax and Sod problems}
\end{figure}

The initial condition are parameterized for $\mu\in [0,1]$ as 
\begin{equation}
\label{euler:3}
V_0(x;\mu)=\mu V_{\text{Sod}}(x)+(1-\mu)V_{\text{Lax}}(x)
\end{equation}
\end{subequations}
and the conservative initial variables $U_0(x;\mu)$ constructed from $V_0(x;\mu)$.

\begin{figure}[h]
\subfigure[Density $\rho$]{\includegraphics[width=0.45\textwidth]{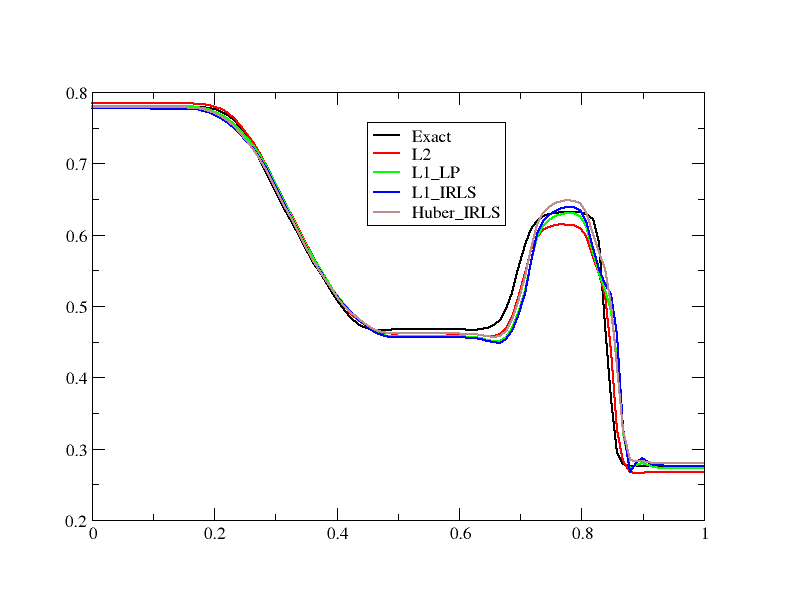}}
 \subfigure[Velocity $u$]{\includegraphics[width=0.45\textwidth]{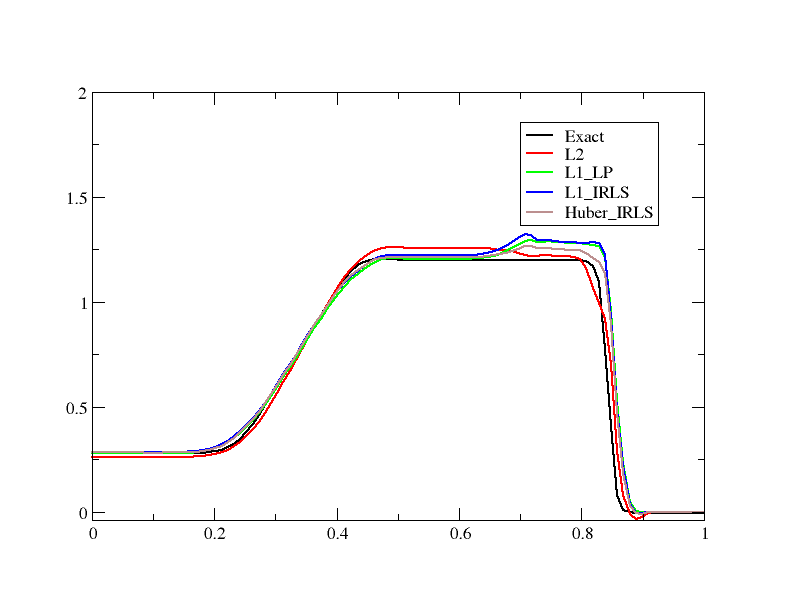} }
 \centering{\subfigure[Pressure $p$]{\includegraphics[width=0.45\textwidth]{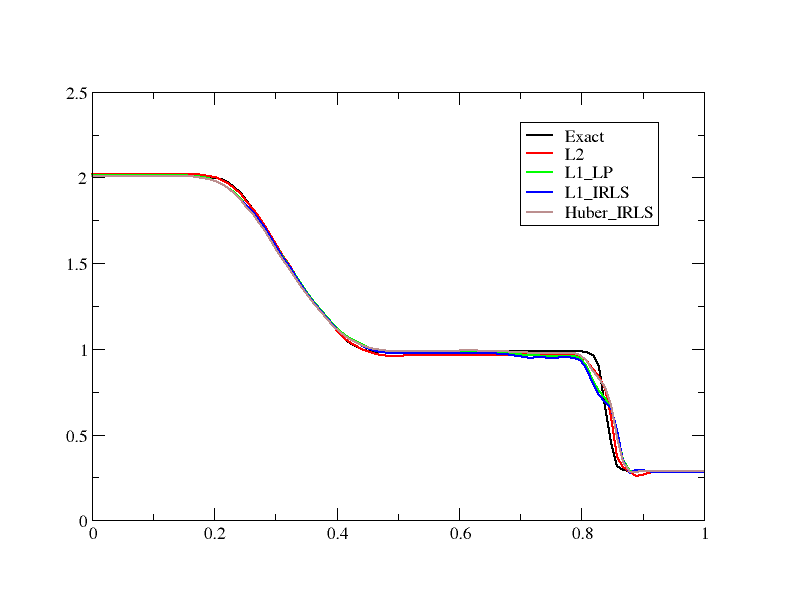} }}
 \caption{\label{fig:reconstruction:un} One-dimensional Euler equations: predicted solutions with strategy  \eqref{reconstruction:un} based on a single expansion}
 \end{figure}

\begin{figure}[h]
\subfigure[Density $\rho$]{\includegraphics[width=0.45\textwidth]{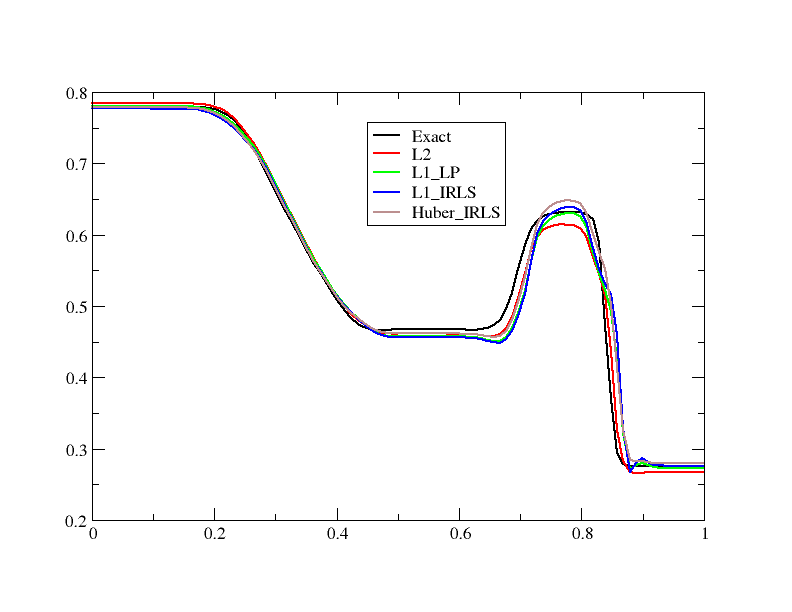}}
 \subfigure[Velocity $u$]{\includegraphics[width=0.45\textwidth]{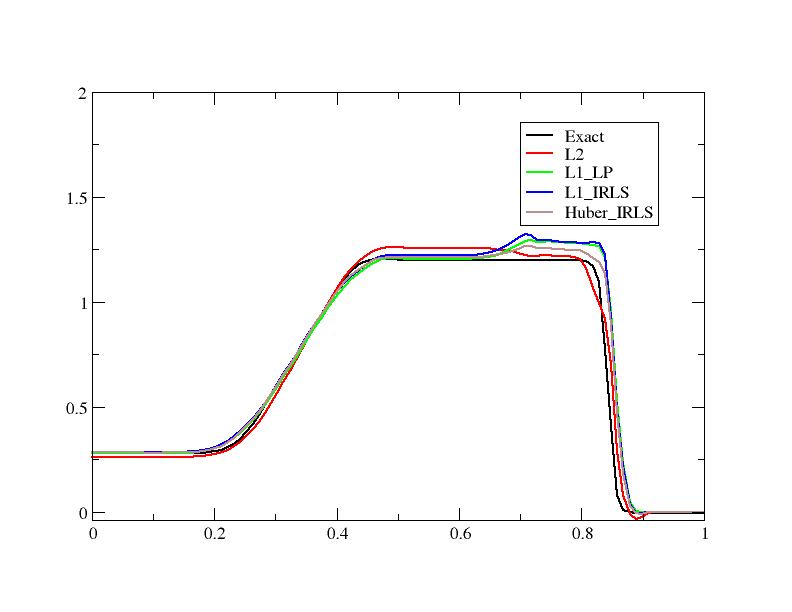} }
 \centering{\subfigure[Pressure $p$]{\includegraphics[width=0.45\textwidth]{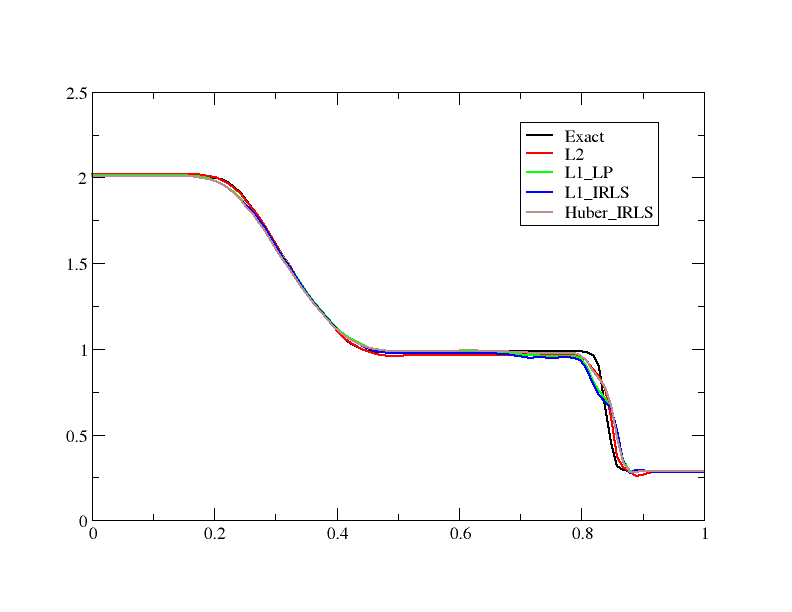} }}
 \caption{\label{fig:reconstruction:trois} One-dimensional Euler equations: predicted solutions with strategy \eqref{reconstruction:trois} based on multiple expansions}
 \end{figure}

In the subsequent numerical experiments, two strategies are exploited to construct, from the dictionary $\mathcal{D}$, the approximation $\ubold^n(\mu)$ of the solution at each time step $n$:
\begin{itemize}
\item Either we reconstruct together the discretized density vectors $\rhobold$, momentum $\mbold=\rhobold \ubold$ and energy $\Ebold$, i.e. the state variable at time $t_n$ using only one coefficient vector $\alphabold^n=(\alpha^n_1,\cdots,\alpha^n_r)$
\begin{equation}
\label{reconstruction:un}
\ubold^n=\begin{pmatrix}
\rhobold^{n}\\
\mbold^n\\
\Ebold^n\end{pmatrix}\approx \sum_{j=1}^r \alpha_j^n  \ubold^n(\mu_j) . 
\end{equation}
Here the $\{\alpha_j^n\}_{j=1}^r$ are obtained by minimizing $J$  on the density components of the state because the density enable to detect fans, contact discontinuities and shocks, contrarily to pressure and velocity which are constant across contact waves. Doing so we expect to control better the numerical oscillations, if any, than with the other physical variables. Similar arguments could be applied with the other conserved variables as well.
\item Alternatively, we reconstruct each conserved variable separately
\begin{equation}\label{reconstruction:trois}
\rhobold^n\approx\sum_{j=1}^N \alpha_{\rho}^n \rhobold^n(\mu_j), \qquad
\mbold^n\approx\sum_{j=1}^N \alpha_{m}^n \mbold^n(\mu_j), \qquad
\Ebold^n\approx\sum_{j=1}^N \alpha_{E}^n \Ebold^n(\mu_j) .
\end{equation}
where the minimization procedures are done \textit{independently} on each conserved variable.
\end{itemize}
In order to test these approaches, the PDE is discretized by finite volumes using a discretization resulting in $Np=3000$ dofs. The parameter range $\mathcal{D}=\{0.0, 0.2,0.4,0.5,0.8,1\}$ is considered together with a target $\mu^\star=0.6$. The results using the first strategy, see eq. \eqref{reconstruction:un}, are displayed in Figure~\ref{fig:reconstruction:un} and those using the second strategy, see eq. \eqref{reconstruction:trois}, reported in Figure \ref{fig:reconstruction:trois}.
 
 From both figures, we can see that the overall structure of the solutions is correct. Nevertheless, there are differences that can be highlighted.
From Figure \ref{fig:reconstruction:un}, we can observe that the density predictions, besides an undershoot at the shock, are good for all minimization procedures except the $L^2$-norm-based one which is oscillatory. However, we cannot recover correct values of the initial velocity (see left  boundary), because there is no reason to believe that the coefficient $\alphabold$, evaluated from the density \emph{only}, will also be correct for the momentum.   A careful observation of the pressure plot also reveals the same behavior which is not satisfactory. For the same reason, if any other \emph{single} variable is used for a global approximation of each conservative variables, there no reason why better qualitative results could be obtained.

This problem does not occur with the second strategy for the reconstruction \eqref{reconstruction:trois}: the correct initial values are recovered. The four minimization procedures behave approximately the same, and we have a slight undershoot/overshoot at the foot of the shock. This does not appear for the $L^2$-norm minimization, but the density is less satisfactory  between the contact and the shock in that case. All methods have slight problems on the velocity, between the contact and the shock. Overall, the Huber function minimization seems to produce the most satisfactory results.

In order to obtain these results we have been faced to the following issue. Take the momentum, for example. For at least half of the mesh points, its value is $0$, and for half of the points, its value is set to a constant. Hence, the matrix $\Abold$ used in the minimization procedure and built on the momentum dictionary has a rank 2 only. The same is true for the other variables, and we are looking here for $r$ coefficients. Two approaches can be followed to address this issue. The first one relies on Gram-Schmidt orthogonalization of the solutions prior to their use as a basis for the solution. The second approach, followed here, consists into perturbing infinitesimally and randomly the matrices involved in the procedure,
so $A_{ij}$ is replaced by $A_{ij}+\varepsilon_{ij}$. The distribution of $\varepsilon_{ij}$ is uniform. This has the effect of giving the maximum possible rank to the perturbed matrix\footnote{We have had to use the same technique for  Burgers' equation in section \ref{Unsteady Burgers}.}. We have expressed that $\epsilon_{ij}$ should depend on the variable, we have chosen
$$\varepsilon_{ij}=\epsilon_{ij} L_{\text{ref}}$$ where $L_{\text{ref}}$ is the difference between the minimum and the maximum, over the dictionary, of the considered variable. Choosing the same $\varepsilon_{ij}$ is taken for all variables, this has the effect of increasing the amplitude of the oscillations after then shock, and our experience is that the $L^2$-norm procedure is more sensitive to the choice of $L_{\text{ref}}$.

All this being said, the solution using three distinct coefficients obtained independently is of significantly much better quality than the one using only one expansion.


\section{Conclusion}\label{sec:conclu}
A novel model reduction that relies on a dictionary approach is developed and tested on several steady and unsteady hyperbolic problems. All of the solutions of the problem tested are parameterized and have regions of their spatial domain with discontinuities, leading to solutions with very distinct behaviors, such as different wave speeds and shock locations, making them challenging to reduce using classical projection-based model reduction techniques. 

To address this challenge, the proposed approach based on a dictionary of solutions is coupled with a functional minimization. The analysis and numerical experiments conducted in this work show that the proposed approach is robust (at least for one-dimensional problems) and performs best when the functional is  of $L^1$-norm-type. Among those, the Huber functional is found to be the most robust in all test cases. For all the functional considered, effective numerical algorithms for their minimization are considered, both in the linear and nonlinear cases. Algorithms based on convex programming are the most computationally expensive ones and, as a result, approaches based on iteratively solving $L^2$-norm minimization problems are considered. These approaches are found to  be inexpensive and lead to similar solutions as convex programming. Furthermore, hyper-reduction approaches developed for the solution of $L^2$-norm minimization~\cite{carlberg11,GNAT,amsallem14:morepas} can be readily applied for their efficient solution. This will be the subject of further work.
Another extension of the present work will be the development of appropriate parameter sampling techniques together with error estimators  to select the dictionary elements. 

For unsteady systems, in the present paper, only dictionary members computed at the same time instant are considered in the reduced approximation.  Variants of the proposed framework will be explored. In particular dictionary members associated with multiple time instants will be considered as well, resulting in local ROM approaches~\cite{dihlmann11,amsallem12:localROB,amsallem14:morepas}.

%
%
%
%
%
%


\section*{Acknowledgements} 
The first author has been funded in part by the MECASIF project  (2013-2017) funded by the French "Fonds Unique Interministériel" and SNF grant \# 200021\_153604 of the Swiss National Foundation.
The last author would like to acknowledge partial support by the Army Research Laboratory through the Army High Performance Computing Research Center under Cooperative Agreement W911NF- 07-2-0027, and partial support by the Office of Naval Research under
grant no. N00014-11-1-0707. This document does not necessarily reflect the position of these institutions, and no official endorsement
should be inferred.

\bibliographystyle{unsrt}
\bibliography{reference}

\appendix
\section{{ Algorithms}}\label{algo}
This section reviews  the minimization algorithms that have been used in this study.
\subsection{{$L^2$-norm minimization}}
\subsubsection{{Linear Case}}\label{L2:L}
{The least-squares problem is defined for a skinny matrix $\Abold=\Mbold\Dbold\in\mathbb{R}^{N\times r}$ and a vector $\bbold\in\mathbb{R}^N$ as
\begin{equation}
\min_{\zbold} \| \Abold\zbold +\bbold\|_2^2.
\end{equation}
One possible solution is by QR decomposition as described in Algorithm~\ref{alg:linL2}.
}
\begin{algorithm}
\fontsize{10pt}{10pt}\selectfont
\caption{Linear $L^2$-norm minimization (Least-squares) by QR decomposition} 
\begin{algorithmic}[1]
{\REQUIRE Matrix $\Abold$ and vector $\bbold$
\ENSURE Solution $\zbold$
 \STATE Compute the QR decomposition of $\Abold$
 \[ \Abold = \Qbold \Rbold\]
 \STATE Let $\zbold = - \Qbold^T \Rbold^{\dagger} \bbold$
 }
 \end{algorithmic}\label{alg:linL2}
\end{algorithm}

\subsubsection{{Nonlinear Case}} \label{L2:NL}
{This algorithm is used in the case drawn in Remark \ref{remark:label2}.
The following nonlinear least-squares problem can be solved by the Gauss-Newton method as described in Algorithm~\ref{alg:nonlinL2}
\begin{equation}
\min_{\zbold} \| \rbold(\zbold)\|_2^2.
\end{equation}
Another approach to solve this problem is by the Levenberg-Marquardt procedure~\cite{nocedalbook}.}
\begin{algorithm}
\fontsize{10pt}{10pt}\selectfont
\caption{Nonlinear $L^2$-norm minimization  by the  Gauss-Newton  method} 
\begin{algorithmic}[1]
{\REQUIRE Residual function $\rbold(\cdot)$ and associated Jacobian $\Jbold(\cdot)$, dictionary $\Dbold$, initial guess $\zbold^0$, tolerance for convergence $\epsilon$
\ENSURE Solution $\zbold$
\STATE $l=0$
\STATE Compute $\rbold^0 =\rbold(\Dbold\zbold^0)$ and $\Wbold^0 =\Jbold(\Dbold\zbold^0)\Dbold$
 \WHILE {$ \left\| \left(\Wbold^l\right)^T\rbold^l\right\|_2 > \epsilon\left\| \left(\Wbold^0\right)^T\rbold^0\right\|_2$} 
 \STATE Solve the linear $L^2$-norm minimization problem using Algorithm~\ref{alg:linL2} with arguments $\Wbold^l$ and $\rbold^l$
 \[ \Delta\zbold^l =  \argmin_{\ybold}\left\| \Wbold^l \ybold + \rbold^l\right\|_2^2 \]
 \STATE $\zbold^{l+1} = \zbold^l + \Delta \zbold^l$
 \STATE Compute $\rbold^{l+1} =\rbold(\Dbold\zbold^{l+1})$ and $\Wbold^{l+1} =\Jbold(\Dbold\zbold^{l+1})\Dbold$
   \STATE $l=l+1$
 \ENDWHILE
 \STATE $\zbold = \zbold^l$
 }
 \end{algorithmic}\label{alg:nonlinL2}
\end{algorithm}

\subsection{$L^1$-norm minimization}
\subsubsection{Linear Case}\label{app:L1lin}
{The linear $L^1$-norm minimization  problem is defined for a skinny matrix $\Abold=\Mbold\Dbold\in\mathbb{R}^{N\times r}$ and a vector $\bbold\in\mathbb{R}^N$ as
\begin{equation}\label{eq:L1lin}
\min_{\zbold} \| \Abold\zbold +\bbold\|_1.
\end{equation}
This problem can be recast as a Linear Program as described in Algorithm~\ref{alg:linL1LP}.}

\begin{algorithm}
\fontsize{10pt}{10pt}\selectfont
\caption{Linear $L^1$-norm minimization by Linear Programming} 
\begin{algorithmic}[1]
{\REQUIRE Matrix $\Abold$ and vector $\bbold$
\ENSURE Solution $\zbold$
 \STATE Solve the linear program
\begin{equation*}\label{eq:LP}
\begin{aligned}
(\zbold^\star,\sbold^\star,\tbold^\star) = \argmin_{\zbold,\sbold,\tbold}~ \boldsymbol{1}^T(\sbold+\tbold)&\\
\text{s.t.}~\Abold\zbold -\sbold + \tbold &=    \bbold\\
            \sbold &\geq \boldsymbol{0}\\
            \tbold &\geq \boldsymbol{0}
\end{aligned}
\end{equation*}
 \STATE Let $\zbold = \zbold^\star$
 }
 \end{algorithmic}\label{alg:linL1LP}
\end{algorithm}

{ An issue associated with the solution of \eqref{eq:L1lin} is the fact that there are $2N+r$ variables and $3N$ constraints, among which $N$ are equality constraints. $N$ is the number of degrees of freedom in the high-dimensional problem and can be very large for fine discretization problems, rendering the linear program solution intractable. 

Alternatively, the $L^1$-norm minimization problem~(\ref{eq:L1lin}) can be solved by \emph{Iteratively Reweighted Least-Squares} (IRLS)~\cite{daubechies08}. This approach proceeds iteratively by solving a sequence of weighted least-squares problem. An advantage of this approach is that its implementation can rely entirely on existing least-squares solvers. Furthermore, its complexity is similar to that of the $L^2$-norm minimization problem. The procedure is presented in Algorithm~\ref{alg:linL1IRLS}.
}

\begin{algorithm}
\fontsize{10pt}{10pt}\selectfont
\caption{Linear $L^1$-norm minimization by Iteratively Reweighted Least-Squares (IRLS)} 
\begin{algorithmic}[1]
{\REQUIRE Matrix $\Abold$ and vector $\bbold$, initial guess $\zbold^0$
\ENSURE Solution $\zbold$
\STATE $l=0$
\STATE Compute $\rbold^0 =\Abold\Dbold\zbold^0+\bbold$ and $\Wbold^0 =\Abold\Dbold$
 \WHILE { $l=0$ OR $ \| \Delta \zbold^{l-1}\|_1 > \epsilon (1+ \| \zbold^{l-1}\|_1)$} 
 \STATE Compute the weights $\Zbold^l = \text{diag}\left({|r^l_i|}^{-\frac{1}{2}}\right)$
 \STATE Solve the linear $L^2$ minimization problem using Algorithm~\ref{alg:linL2} with arguments $\Zbold^l\Wbold^l$ and $\Zbold^l\rbold^l$
 \[ \Delta\zbold^l =  \argmin_{\ybold}\| \Zbold^l\Wbold^l \ybold + \Zbold^l\rbold^l\|_2^2 \]
 \STATE $\zbold^{l+1} = \zbold^l + \Delta \zbold^l$
 \STATE Compute $\rbold^{l+1} =\Abold\Dbold\zbold^{l+1}+\bbold$ and $\Wbold^{l+1} =\Abold\Dbold$
   \STATE $l=l+1$
 \ENDWHILE
 \STATE $\zbold = \zbold^l$
 }
 \end{algorithmic}\label{alg:linL1IRLS}
\end{algorithm}

\subsubsection{Nonlinear Case}\label{app:L1nonlin}
{
The following nonlinear $L^1$-norm minimization problem can be solved by a Gauss-Newton-like procedure. 
\begin{equation}
\min_{\zbold} \| \rbold(\zbold)\|_1.
\end{equation}
The approach relies on the solution of a sequence of linear $L^1$-norm minimization problems. Algorithm~\ref{alg:nonlinL1LP} describes the approach when it relies on Linear Programming and Algorithm~\ref{alg:nonlinL1IRLS} when it relies on IRLS.
}
\begin{algorithm}
\fontsize{10pt}{10pt}\selectfont
\caption{Nonlinear $L^1$-norm minimization  by the  Gauss-Newton  method with LP} 
\begin{algorithmic}[1]
{\REQUIRE Residual function $\rbold(\cdot)$ and associated Jacobian $\Jbold(\cdot)$, dictionary $\Dbold$, initial guess $\zbold^0$, tolerance for convergence $\epsilon$
\ENSURE Solution $\zbold$
\STATE $l=0$
\STATE Compute $\rbold^0 =\rbold(\Dbold\zbold^0)$ and $\Wbold^0 =\Jbold(\Dbold\zbold^0)\Dbold$
 \WHILE {$l=0$ \OR $ \left| \| \Wbold^l \Delta\zbold^l + \rbold^l\|_1 - \| \rbold^l\|_1\right| > \epsilon\| \rbold^0\|_1$} 
 \STATE Solve the linear $L^1$-norm minimization problem using Algorithm~\ref{alg:linL1LP} with arguments $\Wbold^l$ and $\rbold^l$
 \[ \Delta\zbold^l =  \argmin_{\ybold}\| \Wbold^l \ybold + \rbold^l\|_1 \]
 \STATE $\zbold^{l+1} = \zbold^l + \Delta \zbold^l$
 \STATE Compute $\rbold^{l+1} =\rbold(\Dbold\zbold^{l+1})$ and $\Wbold^{l+1} =\Jbold(\Dbold\zbold^{l+1})\Dbold$
   \STATE $l=l+1$
 \ENDWHILE
 \STATE $\zbold = \zbold^l$
}
 \end{algorithmic}\label{alg:nonlinL1LP}
\end{algorithm}

\begin{algorithm}
\fontsize{10pt}{10pt}\selectfont
\caption{Nonlinear $L^1$-norm minimization  by the  Gauss-Newton  method with Iteratively Reweighted Least-Squares} 
\begin{algorithmic}[1]
{\REQUIRE Residual function $\rbold(\cdot)$ and associated Jacobian $\Jbold(\cdot)$, dictionary $\Dbold$, initial guess $\zbold^0$, tolerance for convergence $\epsilon$
\ENSURE Solution $\zbold$
\STATE $l=0$
\STATE Compute $\rbold^0 =\rbold(\Dbold\zbold^0)$ and $\Wbold^0 =\Jbold(\Dbold\zbold^0)\Dbold$
 \WHILE { $l=0$ OR $ \| \Delta \zbold^{l-1}|_1 >  \epsilon (1+ \| \zbold^{l-1}\|_1)$} 
 \STATE Compute the weights  $\Zbold^l = \text{diag}\left({|r^l_i|}^{-\frac{1}{2}}\right)$
 \STATE Solve the linear $L^2$-norm minimization problem using Algorithm~\ref{alg:linL2} with arguments $\Zbold^l\Wbold^l$ and $\Zbold^l\rbold^l$
 \[ \Delta\zbold^l =  \argmin_{\ybold}\| \Zbold^l\Wbold^l \ybold + \Zbold^l\rbold^l\|_2^2 \]
 \STATE $\zbold^{l+1} = \zbold^l + \Delta \zbold^l$
 \STATE Compute $\rbold^{l+1} =\rbold(\Dbold\zbold^{l+1})$ and $\Wbold^{l+1} =\Jbold(\Dbold\zbold^{l+1})\Dbold$
   \STATE $l=l+1$
 \ENDWHILE
 \STATE $\zbold = \zbold^l$
 }
 \end{algorithmic}\label{alg:nonlinL1IRLS}
\end{algorithm}

\subsection{{ Huber function minimization}}\label{ssec:Huber}
{ An issue with $L^1$-norm minimization is the fact that the function $ \zbold \mapsto \| \zbold \|_1$ is non differentiable at $\zbold=0$, causing potential difficulties in the  numerical solution of the minimizer, as shown in the numerical results of Section~\ref{ssec:advec2D}. Alternatively, the minimization of the Huber function can be used. This function behaves similarly to $x^2$ for small values of $x$ and as $|x|$ for large values of $x$. It is also differentiable everywhere.  The Huber function $\phi_M$ is defined as:
\begin{equation}
 \phi_M(x) = \left\{  
 \begin{array}{l l }
 x^2 & \text{if}~ |x|\leq M \\
M(|x|-M) & \text{otherwise} ,
\end{array}
 \right.
\end{equation}

The $L^1$-norm minimization problem is then replaced by
\begin{equation}\label{eq:Huber}
\min_{\zbold} \sum_{i=1}^N \phi_M(r_i(\zbold))
\end{equation}
This problem can also be solved by the IRLS approach as described in Algorithm~\ref{alg:nonlinHuberIRLS}. This approach requires choosing an appropriate value for $M$. In the present work, the following choice has been found to be robust across all applications
$$M = \epsilon_2\max(1,\max(|r_i|))$$
with $\epsilon_2=10^{-6}$.
}
\begin{algorithm}
\fontsize{10pt}{10pt}\selectfont
\caption{Nonlinear Huber function minimization  by Iteratively Reweighted Least-Squares} 
\begin{algorithmic}[1]
{\REQUIRE Residual function $\rbold(\cdot)$ and associated Jacobian $\Jbold(\cdot)$, dictionary $\Dbold$, initial guess $\zbold^0$, tolerance for convergence $\epsilon$
\ENSURE Solution $\zbold$
\STATE $l=0$
\STATE Compute $\rbold^0 =\rbold(\Dbold\zbold^0)$ and $\Wbold^0 =\Jbold(\Dbold\zbold^0)\Dbold$
 \WHILE { $l=0$ OR $ \| \Delta \zbold^{l-1}|_1 >  \epsilon (1+ \| \zbold^{l-1}\|_1)$} 
 \STATE Compute the weights $\Zbold^l = \text{diag}\left(\delta(r_i< M) + M{|r^l_i|}^{-\frac{1}{2}}\delta(r_i\geq M) \right)$
 \STATE Let $M=\epsilon_2\max(1,\max(|r_i|))$
 \STATE Solve the linear $L^2$-norm minimization problem using Algorithm~\ref{alg:linL2} and arguments $\Zbold^l\Wbold^l$ and $\Zbold^l\rbold^l$
 \[ \Delta\qbold^l =  \argmin_{\ybold}\| \Zbold^l\Wbold^l \ybold + \Zbold^l\rbold^l\|_2^2 \]
 \STATE $\qbold^{l+1} = \qbold^l + \Delta \qbold^l$
 \STATE Compute $\rbold^{l+1} =\rbold(\Dbold\zbold^{l+1})$ and $\Wbold^{l+1} =\Jbold(\Dbold\zbold^{l+1})\Dbold$
   \STATE $l=l+1$
 \ENDWHILE
 \STATE $\zbold = \zbold^l$
 }
 \end{algorithmic}\label{alg:nonlinHuberIRLS}
\end{algorithm}

\begin{figure}
\begin{center}
{\includegraphics[width=0.75\textwidth,clip=]{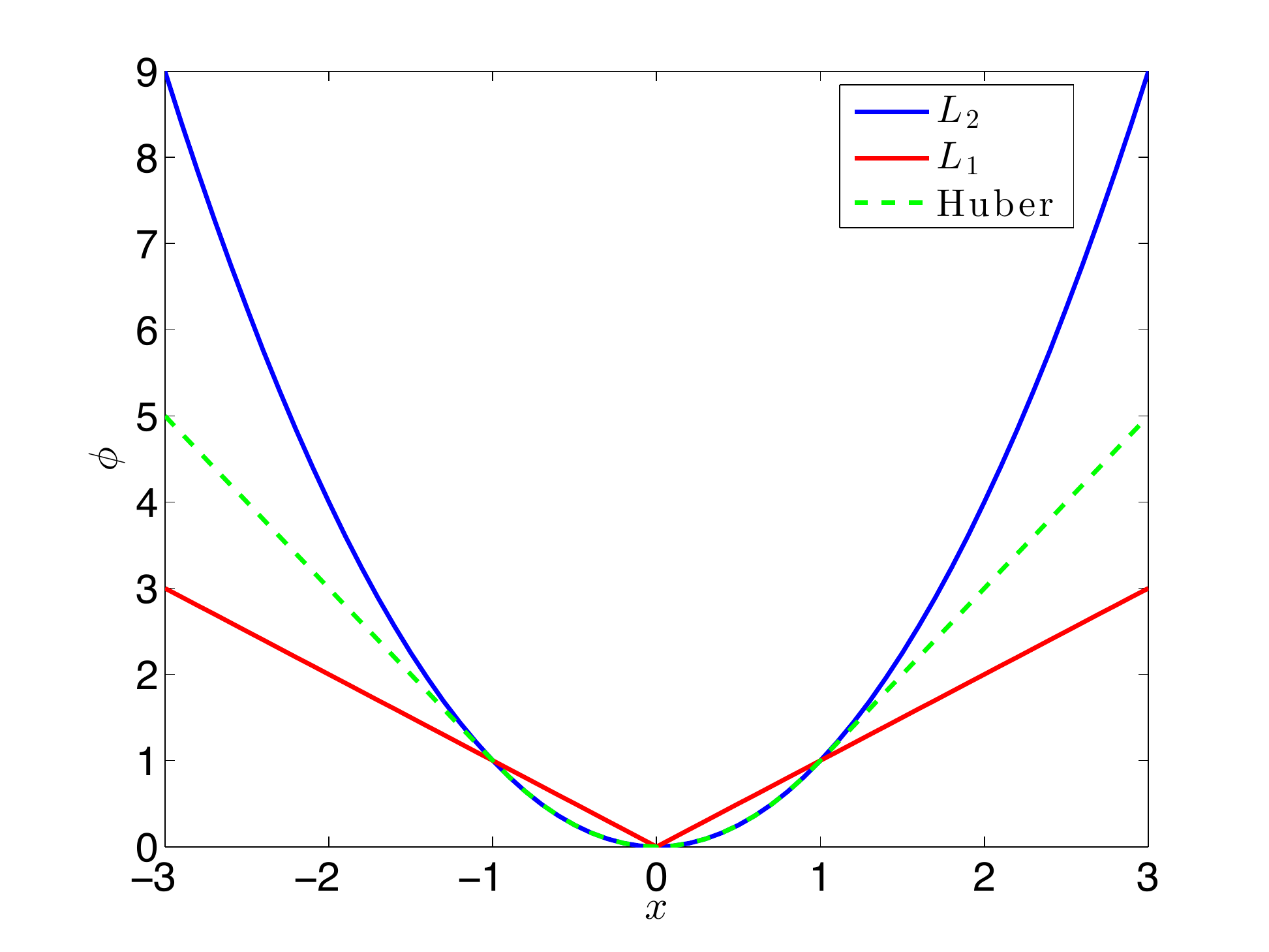}}
\end{center}
\caption{\label{fig:NormsComp}Comparison of the $L^2$, $L^1$ and Huber function norms}
\end{figure}

{ The $L^2$-norm, $L^1$-norm and Huber function minimizations can be all recast with different choices of functions $\phi$ as
\begin{equation}
\min_{\zbold}  \sum_{i=1}^N \phi(r_i(\zbold))
\end{equation}
Figure~\ref{fig:NormsComp} compares the different functions $\phi$ involved.}

\end{document}